\documentclass{amsart}
\usepackage{amssymb}
\usepackage[mathscr]{eucal}

\newtheorem{theorem}{Theorem}[section]

\newtheorem{lemma}[theorem]{Lemma}
\newtheorem{proposition}[theorem]{Proposition}

\theoremstyle{definition}
\newtheorem{definition}[theorem]{Definition}
\newtheorem{notation}[theorem]{Notation}
\newtheorem{example}[theorem]{Example}

\theoremstyle{remark}

\newcommand{\abs}[1]{\lvert#1\rvert}

\begin{document}

\title{Quantized Gromov-Hausdorff distance}

\author{wei wu}
\address{Department of Mathematics, East China Normal University, Shanghai 200062, P.R. China}
\email{wwu@math.ecnu.edu.cn}
\curraddr{Department of Mathematics, University of California, Berkeley, CA 94720-3840}
\email{wwu@math.berkeley.edu}
\subjclass[2000]{Primary 46L87; Secondary 46L07, 53C23, 58B34,
60B10} \keywords{Quantized metric space, matrix Lipschitz
seminorm, matrix seminorm, matrix state space, quantized Gromov-Hausdorff distance}

\begin{abstract}A quantized metric space is a matrix order
unit space equipped with an operator space version of Rieffel's
Lip-norm. We develop for quantized metric spaces an operator
space version of quantum Gromov-Hausdorff distance. We show that
two quantized metric spaces are completely isometric if and
only if their quantized Gromov-Hausdorff distance is zero. We establish a completeness theorem. As applications, we show 
that a quantized metric space with 1-exact underlying matrix order unit space is a limit of matrix algebras with respect to 
quantized Gromov-Hausdorff distance, and that matrix algebras converge naturally to the sphere for quantized Gromov-Hausdorff distance.
\end{abstract}

\maketitle

\section{Introduction}\label{s1}

Following up the compact metric spaces given by Connes in
connection with his theory of quantum Riemannian geometry defined
by Dirac operators \cite{co1}, Rieffel defined the notion of a
compact quantum metric space $(A,L_A)$ in \cite{ri4} as an order
unit space $A$ equipped with a Lip-norm $L_A$, which is a
generalization of the usual Lipschitz seminorm on functions which one associates to an ordinary metric. Many interesting
examples of compact quantum metric space have been constructed
\cite{ri1,ri3,ozri,li1}. Motivated by the type of convergence of
spaces that has recently begun to play a central role in string
theory, Rieffel introduces the quantum Gromov-Hausdorff distance
for the compact quantum metric spaces as a quantum analogue of
Gromov-Hausdorff distance, and shows that the basic theorems of
the classical theory have natural quantum analogues.

In \cite{wu1} and \cite{wu2}, we formulated matrix Lipschitz seminorms on matrix order unit spaces. This operator space version of 
Lipschitz seminorm has many nice properties which are similar to those for ordinary metric spaces. These data may then be thought of as 
some `noncommutative metric spaces'. So
it is natural to ask, as does Rieffel in \cite{ri4}, if it is
possible to develop a corresponding operator space version of
quantum Gromov-Hausdorff distance. This is the aim of the present
article. 

In contrast to the matricial quantum Gromov-Hausdorff distance in \cite{ker} and operator Gromov-Hausdorff distance in \cite{keli}, 
our quantized Gromov-Hausdorff distance operates entirely at the ``matrix" level. Not only the matrix state spaces but also the matrix 
Lipschitz seminorms and the complete isometries are brought into our picture. This should be important in the background of operator 
systems. 

The paper has eight sections. Section 2 contains preliminaries, mainly to fix some terminology and notation. In Section 3 we
define quantized metric space and develop an operator ``quotient". 
Section 4 defines our quantized Gromov-Hausdorff distance, and we prove that it satisfies the
triangle inequality. Section 5 deals with the operator
Gromov-Hausdorff distance zero. We establish that it implies a
complete isometry. Section 6 treats the completeness theorem of the complete isometry classes of quantized metric spaces. 
In Section 7 we show 
that a quantized metric space with 1-exact underlying matrix order unit space is a limit of matrix algebras with respect to 
quantized Gromov-Hausdorff distance. It is established in Section 8 that matrix algebras converge naturally to the sphere for quantized 
Gromov-Hausdorff distance.

\section{Preliminaries}\label{s2}

All vector spaces are assumed to be complex throughout this paper.
Given a vector space $V$, we let $M_{m, n}(V)$ denote the matrix
space of all $m$ by $n$ matrices $v=[v_{ij}]$ with $v_{ij}\in V$,
and we set $M_n(V)=M_{n,n}(V)$. If $V=\mathbb C$,  we write
$M_{m,n}=M_{m,n} (\mathbb C)$ and $M_n=M_{n,n}(\mathbb C)$, which
means that we may identify $M_{m,n}(V)$ with the tensor product
$M_{m,n}\otimes V$. We identify $M_{m,n}$ with the normed space
$\mathcal B({\mathbb C}^n, {\mathbb C}^m)$. We use the standard
matrix multiplication and *-operation for compatible scalar
matrices, and $1_n$ for the identity matrix in $M_n$.

There are two natural operations on the matrix spaces. For $v\in
M_{m,n}(V)$ and $w\in M_{p,q} (V)$, the direct sum $v\oplus w\in
M_{m+p,n+q}(V)$ is defined by letting
\[v\oplus w=\left[\begin{array}{cc}v&0\\ 0&w\end{array}\right],\]
and if we are given $\alpha\in M_{m,p}$, $v\in M_{p,q}(V)$ and $\beta\in M_{q,n}$, the matrix
product  $\alpha v\beta\in M_{m,n}(V)$ is defined by
\[\alpha v\beta=\left[\sum_{k,l}\alpha_{ik}v_{kl}\beta_{lj}\right].\]

A *-{\it vector space} $V$ is a complex vector space together with
a conjugate linear mapping $v\longmapsto v^*$ such that
$v^{**}=v$. A complex vector space $V$ is said to be {\it matrix
ordered} if:
\begin{enumerate}
\item $V$ is a *-vector space;
\item each $M_n(V)$, $n\in\mathbb N$, is partially ordered;
\item  $\gamma^*M_n(V)^+\gamma\subseteq M_m(V)^+$ if $\gamma=[\gamma_{ij}]$ is any $n\times m$
matrix of complex numbers.
\end{enumerate}
A {\it matrix order unit space} $(\mathcal V, 1)$ is a matrix ordered space $\mathcal V$
together with a distinguished order unit $1$ satisfying the following conditions:
\begin{enumerate}
\item $\mathcal V^+$ is a proper cone with the order unit $1$;
\item each of the cones $M_n(\mathcal V)^+$ is Archimedean.
\end{enumerate}
Each matrix order unit space $(\mathcal V, 1)$ may be provided with the norm
\[\|v\|=\inf\left\{t\in\mathbb R:\, \left[\begin{array}{cc} t1&v\\ v^*&t1\end{array}
\right]\ge 0\right\}.\]
As in \cite{ri4}, we will not assume that $\mathcal V$ is complete for the norm.

If $V$ and $W$ are *-vector spaces and $\varphi: V\longmapsto W$
is a linear mapping, we have a linear mapping $\varphi^*:
V\longmapsto W$ defined by $\varphi^*(v)={\varphi(v^*)}^*$.

Given vector spaces $V$ and $W$ and a linear mapping $\varphi:\,
V\longmapsto W$ and $n\in\mathbb N$, we have a corresponding
$\varphi_n:\, M_n(V)\longmapsto M_n(W)$ defined by
\[\varphi_n([v_{ij}])=[\varphi(v_{ij})].\]
If $V$ and $W$ are vector spaces in duality, then they determine
the matrix pairing
\[\ll\cdot, \cdot\gg:\, M_n(V)\times M_m(W)\longmapsto M_{nm},\]
where \[\ll[v_{ij}], [w_{kl}]\gg=\left[<v_{ij}, w_{kl}>\right]\]
for $[v_{ij}]\in M_n(V)$ and $[w_{kl}]\in M_m(W)$.

A {\it graded set} $\mathbf S=(S_n)$ is a sequence of sets
$S_n(n\in\mathbb N)$. If $V$ is a locally convex topological
vector space, then the canonical topology on $M_n(V)(n\in\mathbb
N)$ is that determined  by the natural linear isomorphism
$M_n(V)\cong V^{n^2}$, that is, the product topology. A graded set
$\mathbf S=(S_n)$ with $S_n\subseteq M_n(V)$ is {\it closed} or {\it compact} if
that is the case for each set $S_n$ in the product topology in
$M_n(V)$. Given a vector space $V$, we say that a graded set
$\mathbf B=(B_n)$ with $B_n\subseteq M_n(V)$ is {\it absolutely
matrix convex} if for all $m, n\in\mathbb N$
\begin{enumerate}
\item $B_m\oplus B_n\subseteq B_{m+n}$;
\item $\alpha B_m\beta\subseteq B_n$ for any contractions $\alpha\in M_{n,m}$ and $\beta\in
M_{m,n}$.
\end{enumerate}
A {\it matrix convex set} in $V$ is a graded set
$\mathbf K=(K_n)$ of subsets $K_n\subseteq M_n(V)$ such that 
\[\sum^k_{i=1}\gamma_i^\ast v_i\gamma_i\in K_n\]
for all $v_i\in K_{n_i}$ and $\gamma_i\in M_{n_i,n}$ for $i=1,2,\cdots,k$ satisfying $\sum^k_{i=1}\gamma_i^\ast\gamma_i= 1_n$.
Let $V$ and $W$ be vector spaces in duality, and let $\mathbf
S=(S_n)$ be a graded set with $S_n\subseteq M_n(V)$. The {\it
absolute operator polar} ${\mathbf
S}^{\circledcirc}=(S_n^{\circledcirc})$ with
$S_n^\circledcirc\subseteq M_n(W)$, is defined by
$S_n^\circledcirc=\{w\in M_n(W):\,\|\ll v, w\gg\|\le 1\hbox{ for
all }v\in S_r, r\in\mathbb N\}$. The {\it matrix polar} ${\mathbf S}^{\pi}=(S_n^{\pi})$ with
$S_n^\pi\subseteq M_n(W)$, is defined by
$S_n^\pi=\{w\in M_n(W):\, \mathrm{Re}\ll v, w\gg\le 1_{r\times n}\hbox{ for
all }v\in S_r, r\in\mathbb N\}$. Given a subset $S\subseteq V$, the {\it absolute polar} of $S$ is defined by 
$S^\circ=\{w\in W: |<v,w>|\le 1\mbox{ for all }v\in S\}$.

A {\it gauge} on a vector space $V$ is a function
$g:\,V\longmapsto [0, +\infty]$ such that
\begin{enumerate}
\item $g(v+w)\le g(v)+g(w)$;
\item $g(\alpha v)\le\abs{\alpha} g(v)$,
\end{enumerate}
for all $v,w\in V$ and $\alpha\in\mathbb C$. We say that a gauge $g$ is a {\it seminorm} on
$V$ if $g(v)<+\infty$ for all $v\in V$. Given an arbitrary vector space $V$, a {\it matrix
gauge} $\mathcal G=(g_n)$ on $V$ is a sequence of gauges
\[g_n:\, M_n(V)\longmapsto [0, +\infty]\]
such that
\begin{enumerate}
\item $g_{m+n}(v\oplus w)=\max\{g_m(v), g_n(w)\}$;
\item $g_n(\alpha v\beta)\le\|\alpha\| g_m(v)\|\beta\|$,
\end{enumerate}
for any $v\in M_m(V)$, $w\in M_n(V)$, $\alpha\in M_{n,m}$ and
$\beta\in M_{m,n}$. A matrix gauge $\mathcal G=(g_n)$ is a {\it
matrix seminorm} on $V$ if for any $n\in\mathbb N, g_n(v)<+\infty$
for all $v\in M_n(V)$. If each $g_n$ is a norm on $M_n(V)$, we say
that $\mathcal{G}$ is a {\it matrix norm}. An {\it operator space}
is a vector space together with a matrix norm on it. For a matrix
order unit space $(\mathcal V, 1)$, it is an operator space with
the matrix norm determined by the matrix order on it.

\section{Quantized metric space}\label{s3}

First we recall the following definitions given in \cite{wu1, wu2}:

\begin{definition}\label{def:31} Given a matrix order unit space $(\mathcal V, 1)$, a {\it matrix
Lipschitz seminorm} $\mathcal L$ on $(\mathcal V, 1)$ is a sequence of seminorms
\[L_n:\,M_n(\mathcal V)\longmapsto [0, +\infty)\]
such that
\begin{enumerate}
\item the null space of each $L_n$ is $M_n({\mathbb C1})$;
\item $L_{m+n}(v\oplus w)=\max\{L_m(v), L_n(w)\}$;
\item $L_n(\alpha v\beta)\le\|\alpha\| L_m(v)\|\beta\|$;
\item $L_m(v^*)=L_m(v)$,
\end{enumerate}
for any $v\in M_m(\mathcal V)$, $w\in M_n(\mathcal V)$, $\alpha \in M_{n,m}$ and
$\beta \in M_{m,n}$.
\end{definition}

Let $(\mathcal V, 1)$ be a matrix order unit space. The {\it
matrix state space} of $(\mathcal V, 1)$ is the collection
$\mathcal{CS}(\mathcal V)=(CS_n(\mathcal V))$ of {\it matrix states}
\[CS_n(\mathcal V)=\{\varphi:\,\varphi \hbox{ is a unital completely positive linear mapping from }\mathcal{V}\hbox{ into }M_n\}.\]
If ${\mathcal L}=(L_n)$ is a matrix
Lipschitz seminorm on $(\mathcal V, 1)$, we have a collection $\mathcal{D_L}=(D_{L_n})$ of metrics on 
$\mathcal{CS}(\mathcal V)$ given by
\[D_{L_n}(\varphi, \psi)=\sup\{\|\ll \varphi, a\gg -\ll\psi, a\gg\|:\, a\in M_r(\mathcal V),
L_r(a)\le 1, r\in{\mathbb N}\},\]
for $\varphi , \psi\in CS_n(\mathcal V)$ (notice that it may take value $+\infty$). And in turn we obtain a sequence 
${\mathcal L}_{{\mathcal D}_{\mathcal
L}}=(L_{D_{L_n}})$ of gauges on $(\mathcal V, 1)$ by
\[L_{D_{L_n}}(a)=\sup\left\{\frac{\|\ll \varphi, a\gg -\ll\psi, a\gg\|}{D_{L_r}(\varphi,
\psi)}:\,\varphi, \psi\in CS_r(\mathcal V), \varphi\neq \psi,
r\in{\mathbb N}\right\},\] for all $a\in M_n(\mathcal V)$.

\begin{definition}\label{de31} Let $(\mathcal V, 1)$ be a matrix order unit space.
By a {\it matrix Lip-norm} on $(\mathcal V, 1)$ we mean a matrix
Lipschitz seminorm $\mathcal L=(L_n)$ on $(\mathcal V, 1)$ such
that the $\mathcal{D_L}$-topology on $\mathcal{CS(V)}$ agrees with
the BW-topology.
\end{definition}

We are now prepared to make:

\begin{definition}\label{de32} By a {\it quantized metric space} we mean a pair
$(\mathcal{V}, \mathcal{L})$ consisting of a matrix order unit
space $(\mathcal{V}, 1)$ with a matrix Lip-norm $\mathcal{L}$
defined on it.
\end{definition}

\begin{example}\label{ex33}
Let $(X,\rho)$ be an ordinary compact metric space, let $\mathcal{A}$ denote the set of Lipschitz functions on $X$, and let $L_{\rho}$ 
denote the Lipschitz seminorm on $\mathcal{A}$. Then $\mathcal{A}\subseteq C(X)$, and for $f,g\in\mathcal{A}$ and $\alpha\in\mathbb{C}$, we have
\[L_{\rho}(f^{*})=L_{\rho}(f),\ \ \ \ L_{\rho}(\alpha f)=|\alpha|L_{\rho}(f),\ \ \ \ L_{\rho}(f+g)\le L_{\rho}(f)+L_{\rho}(g).\]
Thus $\mathcal{A}$ is a self-adjoint linear subspace of $C(X)$ which contains constant functions, and so $\mathcal{A}$ is a matrix order unit 
space by Theorem 4.4 in \cite{chef}.

Since $L_{\rho}$ is lower semicontinuous, $K=\{f\in\mathcal{A}: L_{\rho}(f)\le 1\}$ is an absolutely convex normed-closed (and hence 
is weakly closed) set in $\mathcal{A}$. $K$ 
determines a graded set
\[K_n=\left\{\begin{array}{ll}
K,& \mathrm{ if }\ \ n=1,\\
\{0\},& \mathrm{ if }\ \ n>1.\end{array}\right.\]
The minimal envelope $\hat{\mathcal{K}}$ of $K$ is the matrix bipolar $\mathcal{K}^{\circledcirc\circledcirc}$ of $\mathcal{K}$. $\hat{\mathcal{K}}$ 
is an absolutely matrix convex weakly closed graded set. We let $\hat{\mathcal{L}}=(\hat{L}_n)$ be the corresponding matrix gauge of $\hat{\mathcal{K}}$. 
Since $\hat{L}_1=L_\rho$ is a Lipschitz seminorm, $\hat{\mathcal{L}}$ is a matrix Lipschitz seminorm. $\rho_{L_\rho}=\rho$ implies that 
$\hat{\mathcal{L}}$ is also a matrix Lip-norm(see Theorem 1.9 in \cite{ri1} and Proposition 7.5 in \cite{wu2}). Therefore, $(\mathcal{A},\hat{\mathcal{L}})$ is 
a quantized metric space. It is called the {\it minimal quantized metric space} of $(X, \rho)$. The maximal envelope $\check{\mathcal{K}}$ of $K$ is the matrix polar $(\mathcal{K}^{\circ})^\circledcirc$ of 
$\mathcal{K}^\circ=(K^\circ_n)$, where
\[K_n^\circ=\left\{\begin{array}{ll}
K^\circ,& \mathrm{ if }\ \ n=1,\\
\{0\},& \mathrm{ if }\ \ n>1.\end{array}\right.\]
Similarly, $\check{\mathcal{K}}$ is an absolutely matrix convex weakly closed graded set, and the corresponding matrix gauge 
$\check{\mathcal{L}}$ of $\check{\mathcal{K}}$ makes $\mathcal{A}$ into a quantized metric space. $(\mathcal{A},\check{\mathcal{L}})$ is 
called the {\it maximal quantized metric space} of $(X, \rho)$. Moreover, if $\mathcal{C}=(C_n)$ 
is an absolutely matrix convex weakly closed graded set with $C_1=K$, then 
\[\hat{\mathcal{K}}\subseteq\mathcal{C}\subseteq\check{\mathcal{K}},\]
and the corresponding matrix gauge $\mathcal{L}=(L_n)$ satisfies
\[\check{L}_n\le L_n\le\hat{L}_n,\ \ \ \ n\in\mathbb{N},\]
(see page 181 in \cite{efwe}). So $(\mathcal{A}, \mathcal{L})$ is a quantized metric space. It is called a quantized metric space of $(X, \rho)$.
\end{example}

\begin{example}\label{ex34}
Let $(A, L)$ be a compact quantum metric space, that is, an order unit space $(A,e)$ equipped with a seminorm $L$, called Lip-norm, 
on $A$ such that $L(a)=0$ if and only if $a\in\mathbb{R}e$, and the topology on the state space $S(A)$ of $A$ from the metric
\[\rho_L(\mu,\nu)=\sup\{|\mu(a)-\nu(a)|: L(a)\le 1\}\]
is the $w^*$-topology (see Definition 2.2 in \cite{ri4}). So $(S(A), \rho_L)$ is an ordinary compact metric space. 
Let $\mathcal{A}$ denote the set of Lipschitz functions on $S(A)$. By 
Example \ref{ex33}, there exists a quantized metric space structure $(\mathcal{A},\mathcal{L}_1)$ of 
$(S(A),\rho_L)$, where $\mathcal{L}_1=(L_{1,n})$. From Lemma 3.2 in \cite{ri2}, $A\subseteq\mathcal{A}$ and $L_{1,1}(a)\le L(a)$ for $a\in A$. 
Let $\|\cdot\|=(\|\cdot\|_n)$ be the matrix norm determined by the matrix order on $(\mathcal{A},1)$. 
By the basic representation theorem of Kadison\cite{kad}, we also have that $\|a\|=\|a\|_1$ for $a\in A$. If $L$ is lower semicontinuous, the embedding of $A$ into 
$\mathcal{A}$ is {\it isometric}, that is, $\|a\|=\|a\|_1$ and $L(a)=L_{1,1}(a)$ for all $a\in A$, according to Theorem 4.1 in \cite{ri2}.

Set
\[\mathcal{V}=\mathcal{A}\cap(A+iA).\]
We denote the restriction of $\mathcal{L}_1$ on $\mathcal{V}$ by $\mathcal{L}=(L_n)$. Then $\mathcal{V}$ is a self-adjoint linear subspace of $\mathcal{A}$ and contains the order 
unit of $\mathcal{A}$. So $\mathcal{V}$ is a matrix order unit space. Because the $\mathcal{D}_{\mathcal{L}_1}$-topology on $\mathcal{CS}(\mathcal{A})$ 
agrees with the BW-topology, the image of $L_{1,1}^1=\{a\in\mathcal{A}: L_{1,1}(a)\le 1\}$ in $\tilde{\mathcal{A}}=\mathcal{A}/{\mathbb{C}1}$ is totally bounded for $\|\cdot\|^{\sim}_1$ by 
Theorem 5.3 in \cite{wu1}. Since $L_1^1\subseteq L^1_{1,1}$, the image of $L_1^1$ in $\tilde{\mathcal{V}}\subseteq\tilde{\mathcal{A}}$ is totally bounded for 
$\|\cdot\|^{\sim}_1$, and so, by Theorem 5.3 in \cite{wu1}, the $\mathcal{D_L}$-topology on $\mathcal{CS}(\mathcal{V})$ is the BW-topology. Therefore, 
$(\mathcal{V},\mathcal{L})$ is a quantized metric space, and the embedding of $(A,L)$ into $(\mathcal{V},\mathcal{L})$ is an isometry if $L$ is 
lower semicontinuous.
\end{example}

Let $(\mathcal{V}, 1)$ and $(\mathcal{W}, 1)$ be matrix order unit
spaces, and let $\varphi: \mathcal{V}\mapsto\mathcal{W}$ be a
unital completely positive linear mapping. Then we have the dual
mapping $\varphi^{\prime}: \mathcal{W}^*\mapsto\mathcal{V}^*$
determined by $\varphi^{\prime}(f)(v)=f(\varphi(v))$. Let $\tilde{\mathcal{V}}^\ast$ denote the dual space of 
$\tilde{\mathcal{V}}=\mathcal{V}/(\mathbb{C}1)$. $\tilde{\mathcal{V}}^\ast$ is just the subspace of $\mathcal{V}^\ast$ consisting of those 
$f\in\mathcal{V}^\ast$ such that $f(1)=0$. For any $v\in
M_n(\mathcal{V})$ and $g\in M_m(\mathcal{W}^*)$, we have
\[\ll g, \varphi_n(v)\gg=[g_{kl}(\varphi(v_{ij}))]=[(\varphi^{\prime}(g_{kl}))(v_{ij})]=\ll (\varphi^{\prime})_m(g),v\gg.\]
So $\|\varphi\|_{cb}=1$,
$(\varphi^{\prime})_m(M_m(\tilde{\mathcal{W}}^*))\subseteq
M_m(\tilde{\mathcal{V}}^*)$ and
$(\varphi^{\prime})_m(CS_m(\mathcal{W}))\subseteq CS_m(\mathcal{V})$. Moreover, $\varphi^\prime$ is $w^\ast$-continuous.
Let $\varphi_n^{c}=(\varphi^{\prime})_n|_{CS_n(\mathcal{W})}$ for $n\in\mathbb{N}$. Then for $v\in M_n(\mathcal{V})$, $f_i\in
CS_{n_i}(\mathcal{W})$ and $\gamma_i\in M_{n_i,m}$ satisfying
$\sum_{i=1}^{k}\gamma_i^*\gamma_i=1_m$, we have
\[\begin{array}{rcl}
\ll\varphi^c_m\left(\sum_{i=1}^{k}\gamma_i^*f_i\gamma_i\right),v\gg
&=&\ll\sum_{i=1}^{k}\gamma_i^*f_i\gamma_i,\varphi_n(v)\gg\\
&=&\sum_{i=1}^{k}(\gamma_i\otimes 1_n)^*\ll
f_i,\varphi_n(v)\gg(\gamma_i\otimes 1_n)\\
&=&\sum_{i=1}^{k}(\gamma_i\otimes 1_n)^*\ll
\varphi^c_{n_i}(f_i),v\gg(\gamma_i\otimes 1_n)\\
&=&\ll\sum_{i=1}^k\gamma_i^*\varphi^c_{n_i}(f_i)\gamma_i,v\gg.
\end{array}\]
So $\varphi^c=(\varphi_n^c)$ is a BW-continuous matrix affine
mapping of $\mathcal{CS(W)}$ into $\mathcal{CS(V)}$. In
particular,
$\varphi^c(\mathcal{CS(W)})=(\varphi^c_n(CS_n(\mathcal{W})))$ is a
closed matrix convex subset of $\mathcal{CS(V)}$. Clearly
$\varphi^c_n$ is injective if $\varphi$ is surjective.

Let $\mathcal{L}$ be a matrix Lipschitz seminorm on $\mathcal{V}$. On $\mathcal{V}^\ast$, we define 
the matrix gauge $\mathcal{L}^\prime=(L_n^\prime)$ by
\[L_n^\prime(f)=\sup\left\{\|\ll f, a\gg\|: a\in L^1_r, r\in\mathbb{N}\right\},\quad f\in M_n(\mathcal{V}^\ast).\]
Then $L_n^\prime(\varphi-\psi)=D_{L_n}(\varphi,\psi)$ for $\varphi,\psi\in CS_n(\mathcal{V})$ by Lemma 4.3 in \cite{wu2}.

\begin{proposition}\label{pro33} Let $(\mathcal{V}, 1)$ and $(\mathcal{W},
1)$ be matrix order unit spaces, and let $\varphi:
\mathcal{V}\mapsto\mathcal{W}$ be a unital completely positive
linear mapping which is surjective. Let $\mathcal{L}$ be a matrix
Lipschitz seminorm on $\mathcal{V}$, and let
$\mathcal{L_W}=(L_{\mathcal{W},n})$ be a sequence of the
corresponding quotient seminorms on $\mathcal{W}$, defined by
\[L_{\mathcal{W},n}(b)=\inf\{L_n(a): \varphi_n(a)=b\},\ \ \ \ b\in M_n(\mathcal{W}).\]
Then
\begin{enumerate}
\item $\mathcal{L}_{\mathcal{W}}$ is a matrix Lipschitz seminorm on $\mathcal{W}$;
\item $\varphi^{\prime}$ is a complete isometry for the matrix norms
$\mathcal{L}^{\prime}_{\mathcal{W}}$ and $\mathcal{L^{\prime}}$ on
$\tilde{\mathcal{W}}^{*}$ and $\tilde{\mathcal{V}}^{*}$;
\item $\varphi^c$ is a complete isometry for the corresponding
matrix metrics $\mathcal{D_{L_W}}$ and $\mathcal{D_{L}}$;
\item If $\mathcal{L}$ is a matrix Lip-norm, then so is
$\mathcal{L_W}$.
\end{enumerate}
\end{proposition}

\begin{proof} (1) For $b_1\in M_m(\mathcal{W}), b_2\in M_n(\mathcal{W})$, we have
\[\begin{array}{rcl}
&&L_{\mathcal{W},m+n}(b_1\oplus b_2)\\
&=&\inf\{L_{m+n}(a): \varphi_{m+n}(a)=b_1\oplus b_2\}\\
&\le&\inf\{L_{m+n}(a_1\oplus a_2): \varphi_{m+n}(a_1\oplus a_2)=b_1\oplus b_2\}\\
&=&\inf\{\max\{L_m(a_1), L_n(a_2)\}:\varphi_m(a_1)=b_1,
\varphi_n(a_2)=b_2\}\\
&=&\max\{\inf\{L_m(a_1): \varphi_m(a_1)=b_1\},\inf\{L_n(a_2):
\varphi_n(a_2)=b_2\}\}\\
&=&\max\{L_{\mathcal{W},m}(b_1),L_{\mathcal{W},n}(b_2)\}.
\end{array}\]
If $\alpha\in M_{m,n}$, $\beta\in M_{n,m}$ and $b\in
M_n(\mathcal{W})$, we have
\[\begin{array}{rcl}
L_{\mathcal{W},m}(\alpha b\beta) &=&\inf\{L_m(a): \varphi_m(a)=\alpha b\beta\}\\
&\le&\inf\{L_m(\alpha a\beta): \varphi_n(a)=b\}\\
&\le&\|\alpha\|\|\beta\|\inf\{L_n(a): \varphi_n(a)=b\}\\
&=&\|\alpha\|\|\beta\|L_{\mathcal{W},n}(b),
\end{array}\]
and
\[\begin{array}{rcl}
L_{\mathcal{W},n}(b^*)&=&\inf\{L_n(a): \varphi_n(a)=b^*\}\\
&=&\inf\{L_n(a): \varphi_n(a^*)=b\}\\
&=&\inf\{L_n(a^*): \varphi_n(a)=b\}\\
&=&\inf\{L_n(a): \varphi_n(a)=b\}\\
&=&L_{\mathcal{W},n}(b).
\end{array}\]
Given $[\lambda_{ij}]\in M_n$. We have
\[\begin{array}{rcl}
L_{\mathcal{W},n}([\lambda_{ij}1]) &=&\inf\{L_{n}(a): \varphi_n(a)=[\lambda_{ij}1]\}\\
&\le&L_n([\lambda_{ij}1])=0,
\end{array}\]
and so $L_{\mathcal{W},n}([\lambda_{ij}1])=0$. If $b=[b_{ij}]\in M_n(\mathcal{W})$ with $L_{\mathcal{W},n}(b)=0$, then 
\[L_{\mathcal{W},1}(b_{ij})=0,\ \ \ i,j=1,2,\cdots,n.\]
Letting $b_{ij}=c_{ij}+id_{ij}(i,j=1,2,\cdots,n)$, where $c^*_{ij}=c_{ij}, d^*_{ij}=d_{ij}$, we get
\[L_{\mathcal{W},1}(c_{ij})=0,\ \ \ L_{\mathcal{W},1}(d_{ij})=0,\ \ \ i,j=1,2,\cdots,n.\] 
Since $L_{\mathcal{W},1}(c_{ij})=\inf\{L_1(a): \varphi_1(a)=c_{ij}\}$ and $\varphi$ is positive, we have
\[L_{\mathcal{W},1}(c_{ij})=\inf\{L_1(a): \varphi_1(a)=c_{ij}, a=a^*\}.\]
Now by Proposition 3.1 in \cite{ri4}, there exists an $\alpha_{ij}\in\mathbb{R}$ such that $c_{ij}=\alpha_{ij}1$. Similarly, there 
exists a $\beta_{ij}\in\mathbb{R}$ such that $d_{ij}=\beta_{ij}1$. Therefore, $b\in M_n(\mathbb{C}1)$. Thus
$\mathcal{L}_{\mathcal{W}}$ is a matrix Lipschitz seminorm on
$\mathcal{W}$.

(2) Let $f\in M_m(\tilde{\mathcal{W}}^\ast)$. For any $a\in
M_n(\mathcal{V})$ we clearly have
$L_{\mathcal{W},n}(\varphi_n(a))\le L_n(a)$, and so if $L_n(a)\le
1$ we have
\[\|\ll(\varphi^{\prime})_m(f), a\gg\|=\|\ll f, \varphi_n(a)\gg\|\le L^{\prime}_{\mathcal{W},m}(f).\]
Consequently, $L_m^{\prime}((\varphi^{\prime})_m(f))\le
L_{\mathcal{W},m}^{\prime}(f)$. But let $\delta>0$ be given, and let
$b\in M_n(\mathcal{W})$ with $L_{\mathcal{W},m}(b)\le 1$. Then
there is an $a\in M_n(\mathcal{V})$ such that $\varphi_n(a)=b$ and
$L_n(a)\le 1+\delta$. Thus, $L_n(a/(1+\delta))\le 1$.
Consequently,
\[\begin{array}{rcl}
L_m^{\prime}((\varphi^{\prime})_m(f)) &\ge&\|\ll(\varphi^{\prime})_m(f),
a/(1+\delta)\gg\|\\
&=&\|\ll f, \varphi_n(a)\gg\|/(1+\delta)\\
&=&\|\ll f, b\gg\|/(1+\delta).\end{array}\] Taking the supremum
over $b\in M_n(\mathcal{W})$ with $L_{\mathcal{W},n}(b)\le 1$, we
see that
\[L_m^{\prime}((\varphi^{\prime})_m(f))\ge L_{\mathcal{W},m}^{\prime}(f)/(1+\delta).\]
Since $\delta$ is arbitrary, we obtain that
$L_m^{\prime}((\varphi^{\prime})_m(f))\ge L_{\mathcal{W},m}^{\prime}(f)$. Thus
$\varphi^{\prime}$ is a complete isometry.

(3) By (2), we have
\[\begin{array}{rcl}
D_{L_{\mathcal{W},n}}(\phi,
\psi)&=&L^{\prime}_{\mathcal{W},n}(\phi-\psi)\\
&=&L_n^{\prime}((\varphi^{\prime})_n(\phi-\psi))\\
&=&D_{L_n}(\varphi^c_n(\phi),\varphi^c_n(\psi)),\end{array}\]
where $\phi, \psi\in CS_n(\mathcal{W})$, that is, $\varphi^c$ is a
complete isometry for the corresponding matrix metrics
$\mathcal{D_{L_W}}$ and $\mathcal{D_{L}}$.

(4) Suppose that $\mathcal{L}$ is a matrix Lip-norm. Since $\varphi^\prime$ is $w^\ast$-continuous, 
$\varphi$ is surjective, and $CS_n(\mathcal{V})$ is BW-compact, $(\varphi^\prime)_n$ is a homeomorphism of $CS_n(\mathcal{W})$ 
onto $(\varphi^\prime)_n(CS_n(\mathcal{W}))\subseteq CS_n(\mathcal{V})$. Because $D_{L_n}$ gives the BW-topology on $CS_n(\mathcal{V})$, 
$\left.D_{L_n}\right|_{(\varphi^\prime)_n(CS_n(\mathcal{W}))}$ gives the relative topology of $(\varphi^\prime)_n(CS_n(\mathcal{W}))$. 
According to (3), $D_{L_{\mathcal{W},n}}$ gives the BW-topology on $CS_n(\mathcal{W})$. Therefore, $\mathcal{L}_{\mathcal{W}}$ is a matrix Lip-norm.
\end{proof}

\begin{notation}\label{no34} Under the conditions of Propositions \ref{pro33}
we will say that $\mathcal{L}$ induces $\mathcal{L_W}$ via
$\varphi$.
\end{notation}

For a matrix convex set $\mathbf{K}$ in a locally convex vector space, let $A(\mathbf{K})$ be the set of all matrix affine mappings from 
$\mathbf{K}$ to $\mathbb{C}$ (see \S 6 in \cite{wu2}). On the other hand, we have

\begin{proposition}\label{pro335}
Let $(\mathcal{V},1)$ be a matrix order unit space, and let $\mathcal{K}=(K_n)$ be a compact matrix convex subset of $\mathcal{CS}(\mathcal{V})$. 
View the elements of $\mathcal{V}$ as matrix affine mapping from $\mathcal{CS}(\mathcal{V})$ to $\mathbb{C}$ (Proposition 6.1 in \cite{wu2}), and let 
$\mathcal{W}$ consist of their restrictions to $\mathcal{K}$, with $\phi$ the restriction mapping of $\mathcal{V}$ onto $\mathcal{W}$. Then 
$(\mathcal{W}, \phi(1))$ is a matrix order unit space, and $\mathcal{K}=\phi^c(\mathcal{CS}(\mathcal{W}))$.
\end{proposition}

\begin{proof}
Clearly, with the natural matrix order structure on $\mathcal{W}$ and the order unit $\phi(1)$, $(\mathcal{W}, \phi(1))$ is 
a matrix order unit space.

For $\varphi\in K_n\subseteq CS_n(\mathcal{V})$, we define the mapping $\psi: \mathcal{W}\mapsto M_n$ by $\psi(\phi(v))=\varphi(v)$. 
Then $\psi\in CS_n(\mathcal{W})$ and $(\phi^c_n(\psi))(v)=\psi(\phi(v))=\varphi(v)$ for $v\in\mathcal{V}$. Thus 
$K_n\subseteq\phi_n^c(CS_n(\mathcal{W}))$.

Suppose that $\varphi_0\in CS_n(\mathcal{V})$ and $\varphi_0\notin K_n$. By Theorem 1.6 in \cite{wewi}, there is a $v=[v_{ij}]\in M_n(\mathcal{V})$ 
and a self-adjoint $\alpha=[\alpha_{ij}]\in M_n$ such that
\[\mathrm{Re}\ll\varphi,v\gg\le\alpha\otimes 1_r,\]
for all $r\in\mathbb{N}$, $\varphi\in K_r$, and 
\[\mathrm{Re}\ll\varphi_0,v\gg\not\leq\alpha\otimes 1_n.\]
So we obtain $\varphi_n(\mathrm{Re}[\alpha_{ij}1-v_{ij}])\ge0$ for all $r\in\mathbb{N}$ and $\varphi\in K_r$. Thus $\phi_n(\mathrm{Re}[\alpha_{ij}1-v_{ij}])\ge 0$ 
in $\mathcal{W}$. If $\varphi_0=\phi^c_n(\psi_0)$ for some $\psi_0\in CS_n(\mathcal{W})$, we would then have that $\mathrm{Re}\ll\varphi_0,v\gg=\mathrm{Re}\ll\psi_0,\phi_n(v)\gg=
\alpha\otimes 1_n-\mathrm{Re}\ll\psi_0,\phi_n([\alpha_{ij}1-v_{ij}])\gg=\alpha\otimes 1_n-\ll\psi_0,\phi_n(\mathrm{Re}[\alpha_{ij}1-v_{ij}])\gg\le\alpha\otimes 1_n$. Thus, 
$\varphi_0\notin\phi^c(CS_n(\mathcal{W}))$. Therefore, $\mathcal{K}=\phi^c(\mathcal{CS}(\mathcal{W}))$.
\end{proof}

\begin{notation} We will call the matrix order unit space  $(\mathcal{W}, \phi(1))$ in the Proposition \ref{pro335} the {\it quotient} of 
$(\mathcal{V},1)$ with respect to $\mathcal{K}$, and will identify $\mathcal{CS}(\mathcal{W})$ with $\mathcal{K}$. 
When $(\mathcal{V},\mathcal{L})$ is a quantized metric space, 
$(\mathcal{W},\mathcal{L}_{\mathcal{W}})$ is a quantized metric space by Proposition \ref{pro33}. 
$(\mathcal{W}, \mathcal{L}_{\mathcal{W}})$ is called the {\it quotient space} of $(\mathcal{V},\mathcal{L})$ with respect to $\mathcal{K}$ 
and $\phi$. \end{notation}

\begin{proposition}\label{pro35}
Let $(\mathcal{V}_1,1)$, $(\mathcal{V}_2,1)$ and
$(\mathcal{V}_3,1)$ be matrix order unit spaces. Suppose that
$\varphi: \mathcal{V}_1\mapsto\mathcal{V}_2$ and $\psi:
\mathcal{V}_2\mapsto\mathcal{V}_3$ are unital completely positive
linear mappings which are surjective. Denote
$\phi=\psi\circ\varphi$. If $\mathcal{L}$ is a matrix Lipschitz
seminorm on $\mathcal{V}_1$, $\mathcal{L}_{\mathcal{V}_2}$ and
$\mathcal{L}_{\mathcal{V}_3}$ are the induced matrix Lipschitz
seminorms of $\mathcal{L}$ via $\varphi$ and $\phi$, respectively,
and $\mathcal{L}_{\mathcal{V}_{23}}$ is the induced matrix
Lipschitz seminorm of $\mathcal{L}_{\mathcal{V}_2}$ via $\psi$,
then
$\mathcal{L}_{\mathcal{V}_{23}}=\mathcal{L}_{\mathcal{V}_3}$.
\end{proposition}

\begin{proof}
This follows by exactly the same argument used for quantum
Gromov-Hausdorff distance in \cite{ri4}.
\end{proof}

\section{Quantized Gromov-Hausdorff distance}\label{s4}

As in the situation of compact quantum metric spaces, we need a corresponding notion of bridge for estimating distance between quantized 
metric spaces.

Let $(\mathcal{V}_1, \mathcal{L}_1)$ and $(\mathcal{V}_2, \mathcal{L}_2)$ be two quantized metric spaces with the matrix norms 
$\|\cdot\|_1=(\|\cdot\|_{1,n})$ and $\|\cdot\|_2=(\|\cdot\|_{2,n})$ determined by their matrix orders on $(\mathcal{V}_1,1)$
and $(\mathcal{V}_2,1)$, respectively. We form the direct sum $\mathcal{V}_1\oplus\mathcal{V}_2$ of operator spaces (see \S 2.6 in \cite{pi1}). 
$(\mathcal{V}_1\oplus\mathcal{V}_2, (1, 1))$ becomes a matrix order unit space. 

\begin{definition}\label{de43}
Let $(\mathcal{V}_1,\mathcal{L}_1)$ and
$(\mathcal{V}_2,\mathcal{L}_2)$ be quantized metric spaces.
A {\it matrix bridge} between $(\mathcal{V}_1,\mathcal{L}_1)$ and
$(\mathcal{V}_2,\mathcal{L}_2)$ is a matrix seminorm $\mathcal{N}$
on $\mathcal{V}_1\oplus\mathcal{V}_2$ such that
\begin{enumerate}
\item $\mathcal{N}$ is matrix continuous for the matrix norm $\|\cdot\|$ on
$\mathcal{V}_1\oplus\mathcal{V}_2$, that is, each $N_n$ is
continuous for $\|\cdot\|_n$ on
$M_n(\mathcal{V}_1\oplus\mathcal{V}_2)$.
\item $N_n((a,b)^\ast)=N_n(a,b)$ for $a\in M_n(\mathcal{V}_1)$ and $b\in M_n(\mathcal{V}_2)$ and $n\in\mathbb{N}$.
\item $N_1(1,1)=0$ but $N_1(1,0)\neq 0$.
\item For any $n\in\mathbb{N}$, $a\in M_n(\mathcal{V}_1)$ and $\epsilon>0$,
there is a $b\in M_n(\mathcal{V}_2)$ such that
\[\max\{L_{2,n}(b),N_n(a,b)\}\le L_{1,n}(a)+\epsilon,\]
and similarly for $\mathcal{V}_1$ and $\mathcal{V}_2$
interchanged.
\end{enumerate}
\end{definition}

\begin{example}\label{ex44}
Suppose $(\mathcal{V}_1,\mathcal{L}_1)$ and
$(\mathcal{V}_2,\mathcal{L}_2)$ are quantized metric spaces. 
Choose $\varphi_1\in CS_1(\mathcal{V}_1)$ and $\psi_1\in CS_1(\mathcal{V}_2)$. For 
$n\in\mathbb{N}$, we define $N_n: M_n(\mathcal{V}_1\oplus\mathcal{V}_2)\mapsto [0,+\infty)$ by
\[N_n(a,b)=\|\ll\varphi_1,a\gg-\ll\psi_1,b\gg\|.\]
Then $\mathcal{N}=(N_n)$ is a matrix seminorm on $\mathcal{V}_1\oplus\mathcal{V}_2$, and satisfies 
the conditions (1), (2) and (3) of Definition \ref{de43}. For any $a\in M_n(\mathcal{V}_1)$ and $\epsilon>0$, 
choose $b=[\varphi_1(a_{ij})1]\in M_n(\mathcal{V}_2)$. Then we have 
\[\max\{L_{2,n}(b), N_n(a,b)\}=0\le L_{1,n}(a)+\epsilon,\]
and similarly if we are given $b\in M_n(\mathcal{V}_2)$. So $\mathcal{N}$ is a matrix bridge between 
$(\mathcal{V}_1,\mathcal{L}_1)$ and $(\mathcal{V}_2,\mathcal{L}_2)$.
\end{example}

\begin{proposition}\label{pro52}
Let $\mathcal{N}$ be a matrix bridge between quantized
metric spaces $(\mathcal{V}_1,\mathcal{L}_1)$ and
$(\mathcal{V}_2,\mathcal{L}_2)$. Define $\mathcal{L}=(L_n)$ on
$\mathcal{V}_1\oplus\mathcal{V}_2$ by
\[L_n(a,b)=\max\{L_{1,n}(a), L_{2,n}(b), N_{n}(a,b)\},\ \ \ a\in M_n(\mathcal{V}_1), b\in M_n(\mathcal{V}_2), n\in\mathbb{N}.\]
Let $\pi_1$ and $\pi_2$ be the projections from $\mathcal{V}_1\oplus\mathcal{V}_2$ onto $\mathcal{V}_1$ and 
$\mathcal{V}_2$, respectively, which are unital completely positive linear surjective mappings.
Then $\mathcal{L}$ is a matrix Lip-norm on $(\mathcal{V}_1\oplus\mathcal{V}_2, (1, 1))$, and it induces $\mathcal{V}_1$ and $\mathcal{V}_2$ 
via $\pi_1$ and $\pi_2$, respectively. If $\mathcal{L}_1$ and $\mathcal{L}_2$ are lower semicontinuous, then
so is $\mathcal{L}$.
\end{proposition}

\begin{proof}
For $a_i\in M_n(\mathcal{V}_i)$ and $b_i\in M_m(\mathcal{V}_i), i=1,2$, we have
\[\begin{array}{rcl}
&&L_{n+m}(a_1\oplus b_1, a_2\oplus b_2)\\
&=&\max\{L_{1,n+m}(a_1\oplus b_1), L_{2,n+m}(a_2\oplus b_2), N_{n+m}(a_1\oplus b_1, a_2\oplus b_2)\}\\
&=&\max\{\max\{L_{1,n}(a_1),L_{1,m}(b_1)\}, \max\{L_{2,n}(a_2), L_{2,m}(b_2)\},\\
&&\max\{N_n(a_1,a_2), N_m(b_1,b_2)\}\}\\
&=&\max\{\max\{L_{1,n}(a_1),L_{2,n}(a_2), N_n(a_1,a_2)\},\\
&&\max\{L_{1,m}(b_1),L_{2,m}(b_2),N_m(b_1,b_2)\}\}\\
&=&\max\{L_n(a_1,a_2),L_m(b_1,b_2)\},
\end{array}\]
and
\[\begin{array}{rcl}
L_n((a_1,a_2)^*)&=&L_n(a_1^*,a_2^*)\\
&=&\max\{L_{1,n}(a_1^*),L_{2,n}(a_2^*),N_n(a_1^*,a_2^*)\}\\
&=&\max\{L_{1,n}(a_1),L_{2,n}(a_2),N_n(a_1,a_2)\}\\
&=&L_n(a_1,a_2),
\end{array}\]
and for $\alpha\in M_{m,n}$ and $\beta\in M_{n,m}$, we have
\[\begin{array}{rcl}
&&L_m(\alpha(a_1,a_2)\beta)\\
&=&L_m(\alpha a_1\beta,\alpha a_2\beta)\\
&=&\max\{L_{1,m}(\alpha a_1\beta),L_{2,m}(\alpha a_2\beta),N_m(\alpha a_1\beta,\alpha a_2\beta)\}\\
&\le&\max\{\|\alpha\| L_{1,n}(a_1)\|\beta\|,\|\alpha\|L_{2,n}(a_2)\|\beta\|,\|\alpha\|N_n(a_1,a_2)\|\beta\|\}\\
&=&\|\alpha\|L_n(a_1,a_2)\|\beta\|.
\end{array}\]
Thus $\mathcal{L}$ is a matrix seminorm. Since
\[L_1(a_{st})\le L_n([a_{ij}])\le\sum_{i,j=1}^n L_1(a_{ij})\]
for $s,t=1,2,\cdots,n$ and $[a_{ij}]\in M_n(\mathcal{V}_1\oplus\mathcal{V}_2)$,
$L_n([\lambda_{ij}(1,1)])=0$ for $[\lambda_{ij}]\in M_n$. If
$L_n([(a_{ij},b_{ij})])=0$, then
$L_{1,n}([a_{ij}])=L_{2,n}([b_{ij}])=0$, and hence
$a_{ij}=\lambda_{ij}1$ and $b_{ij}=\mu_{ij}1$, $i,j=1,2,\cdots,n$,
where $\lambda_{ij}, \mu_{ij}\in\mathbb{C}$. From
$N_n([a_{ij},b_{ij}])=0$ and $N_1(a_{st},b_{st})\le N_n([a_{ij},b_{ij}])$ for $s,t=1,2,\cdots,n$, we have
\[N_1(\lambda_{ij}1,\mu_{ij}1)=0,\ \ \ i,j=1,2,\cdots,n,\]
and so for $i,j=1,2,\cdots,n$,
\[\begin{array}{rcl}
|\lambda_{ij}-\mu_{ij}|N_1(1,0)&=&N_1((\lambda_{ij}-\mu_{ij})1,0)\\
&=&N_1((\lambda_{ij}1,\mu_{ij}1)-(\mu_{ij}1,\mu_{ij}1))\\
&\le&N_1(\lambda_{ij}1,\mu_{ij}1)+N_1(\mu_{ij}1,\mu_{ij}1)\\
&=&0.
\end{array}\]
Thus
$[(a_{ij},b_{ij})]=[(\lambda_{ij}1,\lambda_{ij}1)]=[\lambda_{ij}(1,1)]$.
So $\mathcal{L}$ is a matrix Lipschitz seminorm.

Similar to the same argument used in Theorem 5.2 of \cite{ri4}, we have that $\mathcal{L}$ induces $\mathcal{L}_{1}$ and $\mathcal{L}_2$ 
via $\pi_1$ and $\pi_2$, respectively. By Proposition 3.1 in \cite{wu1}, Proposition 7.5 in \cite{wu2} and Theorem 5.2 in \cite{ri4} (see also \S 2 
in \cite{ri4}), the $\mathcal{D_L}$-topology on
$\mathcal{CS}(\mathcal{V}_1\oplus\mathcal{V}_2)$ agrees with the
BW-topology. Therefore, $\mathcal{L}$ is a matrix Lip-norm on $(\mathcal{V}_1\oplus\mathcal{V}_2, (1, 1))$.

Suppose that $\mathcal{L}_1$ and $\mathcal{L}_2$ are lower
semicontinuous. Clearly, $\mathcal{L}$ is lower semicontinuous
since $\mathcal{N}$ is matrix continuous.
\end{proof}

\begin{notation}\label{no342} We will denote by $\mathcal{M}(\mathcal{L}_1,
\mathcal{L}_2)$ the set of matrix Lip-norms on
$\mathcal{V}_1\oplus\mathcal{V}_2$ which induce both
$\mathcal{L}_1$ and $\mathcal{L}_2$ via $\pi_1$ and $\pi_2$,
respectively. By Proposition \ref{pro52} and Example \ref{ex44}, $\mathcal{M}(\mathcal{L}_1, \mathcal{L}_2)\neq\emptyset$. 
From Proposition \ref{pro33}, we can view $\mathcal{CS}(\mathcal{V}_1)$ and
$\mathcal{CS}(\mathcal{V}_2)$ as closed matrix convex subsets of
$\mathcal{CS}(\mathcal{V}_1\oplus \mathcal{V}_2)$.\end{notation}

Now we introduce our notion of distance for quantized metric spaces.

\begin{definition}\label{de41}
Let $(\mathcal{V}_1, \mathcal{L}_1)$ and $(\mathcal{V}_2,
\mathcal{L}_2)$ be quantized metric spaces. We define the
{\it quantized Gromov-Hausdorff distance}
$\mathrm{dist}_{NC}(\mathcal{V}_1, \mathcal{V}_2)$ between them by
\[\mathrm{dist}_{NC}(\mathcal{V}_1, \mathcal{V}_2)=\inf_{\mathcal{L}=
(L_n)\in\mathcal{M}(\mathcal{L}_1,
\mathcal{L}_2)}\sup_{n\in\mathbb{N}}\left\{n^{-2}\mathrm{dist}_H^{D_{L_n}}\left(CS_n(\mathcal{V}_1),
CS_n(\mathcal{V}_2)\right)\right\},\] where
$\mathrm{dist}_H^{D_{L_n}}(CS_n(\mathcal{V}_1),
CS_n(\mathcal{V}_2))$ is the Hausdorff distance between
$CS_n(\mathcal{V}_1)$ and $CS_n(\mathcal{V}_2)$ for $D_{L_n}$.
\end{definition}

Given a quantized metric space $(\mathcal{V},\mathcal{L})$, we define its {\it diameter} $\mathrm{diam}(\mathcal{V},\mathcal{L})$ to be the diameter of 
$CS_1(\mathcal{V})$ with respect to $D_{L_1}$. The following proposition indicates that the quantized Gromov-Hausdorff distance is always 
finite.

\begin{proposition}\label{pro44}
Let $(\mathcal{V}_1,\mathcal{L}_1)$ and $(\mathcal{V}_2,\mathcal{L}_2)$ be quantized metric spaces. Then
$$\mathrm{dist}_{NC}(\mathcal{V}_1,\mathcal{V}_2)\le2(\mathrm{diam}(\mathcal{V}_1,\mathcal{L}_1)+
\mathrm{diam}(\mathcal{V}_2,\mathcal{L}_2)).$$
\end{proposition}

\begin{proof}
Choosing arbitrarily $\alpha >0$, $\varphi_0\in CS_1(\mathcal{V}_1), \psi_0\in CS_1(\mathcal{V}_2)$, we set
$$N_n(a,b)={\alpha}^{-1}\|\ll\varphi_0,a\gg-\ll\psi_0,b\gg\|,\ \ \ \ a\in M_n(\mathcal{V}_1), 
b\in M_n(\mathcal{V}_2), n\in{\mathbb N}.$$
As Example \ref{ex44}, $\mathcal{N}=(N_n)$ is a matrix bridge between 
$(\mathcal{V}_1,\mathcal{L}_1)$ and $(\mathcal{V}_2,\mathcal{L}_2)$. By Proposition \ref{pro52}, 
$\mathcal{L}=(L_n)$, where
\[L_n(a,b)=\max\{L_{1,n}(a), L_{2,n}(b), N_{n}(a,b)\},\ \ \ a\in M_n(\mathcal{V}_1), b\in M_n(\mathcal{V}_2), n\in\mathbb{N},\]
is a matrix Lip-norm in $\mathcal{M}(\mathcal{L}_1,\mathcal{L}_2)$. Then for $\varphi\in CS_n(\mathcal{V}_1), \psi\in CS_n(\mathcal{V}_2)$, and $(a,b)\in
M_n(\mathcal{V}_1\oplus\mathcal{V}_2)$ with $L_n(a,b)\le 1$, we have
\[\begin{array}{rcl}
&&\|\ll\varphi,a\gg-\ll\psi,b\gg\|\\
&\le&\|\ll\varphi,a\gg-\ll\underbrace{\varphi_0\oplus\cdots\oplus\varphi_0}_n,a\gg\|\\
&&+\|\ll\underbrace{\varphi_0\oplus\cdots\oplus\varphi_0}_n,a\gg-\ll\underbrace{\psi_0\oplus\cdots\oplus\psi_0}_n,b\gg\|\\
&&+\|\ll\underbrace{\psi_0\oplus\cdots\oplus\psi_0}_n,b\gg-\ll\psi,b\gg\|\\
&\le&\sum_{i,j}\|\ll\varphi_{ij}-\delta_{ij}\varphi_{0},a\gg\|+\alpha+\sum_{i,j}\|\ll\psi_{ij}-\delta_{ij}\psi_{0},b\gg\|
\end{array}\]
If $n=1$, we get 
$$\|\ll\varphi,a\gg-\ll\psi,b\gg\|\le\mathrm{diam}(\mathcal{V}_1,\mathcal{L}_1)+\alpha+\mathrm{diam}(\mathcal{V}_2,\mathcal{L}_2),$$
by Proposition 3.1 in \cite{wu1}. If $n>1$, similar to the proof of Proposition 4.2 in \cite{wu1}, there are $\varphi^{(k)}_{ij}\in CS_1(\mathcal{V}_1), k=1,2,3,4$, such that
\[\varphi_{ij}-\delta_{ij}\varphi_{0}=\varphi^{(1)}_{ij}-\varphi^{(2)}_{ij}+i(\varphi^{(3)}_{ij}-\varphi^{(4)}_{ij}).\]
Since $L_{1,n}(a)\le L_n(a,b)\le 1$, we obtain
\[\begin{array}{rcl}
&&\sum_{i,j}\|\ll\varphi_{ij}-\delta_{ij}\varphi_{0},a\gg\|\\
&\le&\sum_{i,j}(\|\ll\varphi^{(1)}_{ij},a\gg-\ll\varphi^{(2)}_{ij},a\gg\|+
\|\ll\varphi^{(3)}_{ij},a\gg-\ll\varphi^{(4)}_{ij},a\gg\|)\\
&\le&\sum_{i,j}(D_{L_{1,1}}(\varphi^{(1)}_{ij},\varphi^{(2)}_{ij})+D_{L_{1,1}}(\varphi^{(3)}_{ij},\varphi^{(4)}_{ij}))\\
&\le&2n^2\mathrm{diam}(\mathcal{V}_1,\mathcal{L}_1).
\end{array}\]  
Applying the same argument, we have
\[\sum_{i,j}\|\ll\psi_{ij}-\delta_{ij}\psi_{0},b\gg\|\le2n^2\mathrm{diam}(\mathcal{V}_2,\mathcal{L}_2).\]
Hence 
\[\|\ll\varphi,a\gg-\ll\psi,b\gg\|\le2n^2(\mathrm{diam}(\mathcal{V}_1,\mathcal{L}_1)+\alpha+\mathrm{diam}(\mathcal{V}_2,\mathcal{L}_2)).\]
The arbitrariness of $\alpha$ implies that 
$\mathrm{dist}_{NC}(\mathcal{V}_1,\mathcal{V}_2)\le2(\mathrm{diam}(\mathcal{V}_1)+\mathrm{diam}(\mathcal{V}_2))$ 
by Proposition 3.1 in \cite{wu1}.
\end{proof}

It is clear that the quantized Gromov-Hausdorff distance is symmetric in $\mathcal{V}_1$ and $\mathcal{V}_2$. We come to prove that it 
satisfies the triangle inequality.

\begin{theorem}\label{th42}
If $(\mathcal{V}_1, \mathcal{L}_1)$, $(\mathcal{V}_2,
\mathcal{L}_2)$ and $(\mathcal{V}_3, \mathcal{L}_3)$ be
quantized metric spaces, then
\[\mathrm{dist}_{NC}(\mathcal{V}_1, \mathcal{V}_3)\le \mathrm{dist}_{NC}(\mathcal{V}_1,
\mathcal{V}_2)+\mathrm{dist}_{NC}(\mathcal{V}_2, \mathcal{V}_3).\]
\end{theorem}

\begin{proof} Given $\epsilon>0$. Then there are
$\mathcal{L}_{12}\in\mathcal{M}(\mathcal{L}_1, \mathcal{L}_2)$ and
$\mathcal{L}_{23}\in\mathcal{M}(\mathcal{L}_2, \mathcal{L}_3)$
such that
\[\sup_{n\in\mathbb{N}}\left\{n^{-2}\mathrm{dist}_H^{D_{L_{12,n}}}\left(CS_n(\mathcal{V}_1),
CS_n(\mathcal{V}_2)\right)\right\}\le
\mathrm{dist}_{NC}(\mathcal{V}_1, \mathcal{V}_2)+\epsilon\] and
\[\sup_{n\in\mathbb{N}}\left\{n^{-2}\mathrm{dist}_H^{D_{L_{23,n}}}\left(CS_n(\mathcal{V}_2),
CS_n(\mathcal{V}_3)\right)\right\}\le
\mathrm{dist}_{NC}(\mathcal{V}_2, \mathcal{V}_3)+\epsilon.\] We
define $\mathcal{L}=(L_n)$ on
$\mathcal{V}_1\oplus\mathcal{V}_2\oplus\mathcal{V}_3$ by
\[L_n(a_1,a_2,a_3)=\max\{L_{12,n}(a_1,a_2), L_{23,n}(a_2,a_3)\}.\]
Then for $a_i\in M_n(\mathcal{V}_i)$ and $b_i\in
M_m(\mathcal{V}_i), i=1,2,3$, we have
\[\begin{array}{rcl}
&&L_{n+m}(a_1\oplus b_1, a_2\oplus b_2, a_3\oplus b_3)\\
&=&\max\{L_{12,n+m}(a_1\oplus b_1,a_2\oplus b_2),
L_{23,n+m}(a_2\oplus b_2,a_3\oplus b_3)\}\\
&=&\max\{\max\{L_{12,n}(a_1,a_2),L_{12,m}(b_1,b_2)\}, \max\{L_{23,n}(a_2,a_3),L_{23,m}(b_2,b_3)\}\}\\
&=&\max\{\max\{L_{12,n}(a_1,a_2),L_{23,n}(a_2,a_3)\}, \max\{L_{12,m}(b_1,b_2),L_{23,m}(b_2,b_3)\}\}\\
&=&\max\{L_n(a_1,a_2,a_3),L_m(b_1,b_2,b_3)\},
\end{array}\]
and
\[\begin{array}{rcl}
L_n((a_1,a_2,a_3)^*)&=&L_n(a_1^*,a_2^*,a_3^*)\\
&=&\max\{L_{12,n}(a_1^*,a_2^*),L_{23,n}(a_2^*,a_3^*)\}\\
&=&\max\{L_{12,n}(a_1,a_2),L_{23,n}(a_2,a_3)\}\\
&=&L_n(a_1,a_2,a_3),
\end{array}\]
and for $\alpha\in M_{m,n}$ and $\beta\in M_{n,m}$, we have
\[\begin{array}{rcl}
L_m(\alpha(a_1,a_2,a_3)\beta) &=&L_m(\alpha a_1\beta,\alpha
a_2\beta,\alpha a_3\beta)\\
&=&\max\{L_{12,m}(\alpha a_1\beta,\alpha a_2\beta),L_{23,m}(\alpha a_2\beta,\alpha a_3\beta)\}\\
&\le&\max\{\|\alpha\|
L_{12,n}(a_1,a_2)\|\beta\|,\|\alpha\|L_{23,n}(a_2,a_3)\|\beta\|\}\\
&=&\|\alpha\|L_n(a_1,a_2,a_3)\|\beta\|.
\end{array}\]
$L_n(a_1,a_2,a_3)=0$ if and only if $L_{12,n}(a_1,a_2)=0$ and
$L_{23,n}(a_2,a_3)=0$, and this is equivalent to that $(a_1,a_2,a_3)\in M_n(\mathbb{C}(1,1,1))$. Therefore, $\mathcal{L}$ is a matrix
Lipschitz seminorm.

Similar to the same argument used in Lemma 4.4 of \cite{ri4}, we have that $\mathcal{L}$ induces $\mathcal{L}_{12}$, 
$\mathcal{L}_{23}$, $\mathcal{L}_{1}$, $\mathcal{L}_{2}$ and $\mathcal{L}_3$ 
for the evident quotient mappings by Proposition \ref{pro35}. By Proposition 3.1 in \cite{wu1}, Proposition 7.5 in \cite{wu2} and Lemma 4.4 in \cite{ri4} (see also \S 2 
in \cite{ri4}), the $\mathcal{D_L}$-topology on
$\mathcal{CS}(\mathcal{V}_1\oplus\mathcal{V}_2\oplus\mathcal{V}_3)$ agrees with the
BW-topology. So $\mathcal{L}$ is a matrix Lip-norm on $(\mathcal{V}_1\oplus\mathcal{V}_2\oplus\mathcal{V}_3, (1, 1, 1))$.

By Proposition \ref{pro33}, we have
\[\sup_{n\in\mathbb{N}}\left\{n^{-2}\mathrm{dist}_H^{D_{L_n}}\left(CS_n(\mathcal{V}_1),
CS_n(\mathcal{V}_2)\right)\right\}\le
\mathrm{dist}_{NC}(\mathcal{V}_1, \mathcal{V}_2)+\epsilon,\]
\[\sup_{n\in\mathbb{N}}\left\{n^{-2}\mathrm{dist}_H^{D_{L_n}}\left(CS_n(\mathcal{V}_2),
CS_n(\mathcal{V}_3)\right)\right\}\le
\mathrm{dist}_{NC}(\mathcal{V}_2, \mathcal{V}_3)+\epsilon\] and
\[\mathrm{dist}_{NC}(\mathcal{V}_1,\mathcal{V}_3)\le\sup_{n\in\mathbb{N}}\left\{n^{-2}\mathrm{dist}_H^{D_{L_n}}(CS_n(\mathcal{V}_1),CS_n(\mathcal{V}_3))\right\}.\]
So
\[\begin{array}{rcl}
&&\mathrm{dist}_{NC}(\mathcal{V}_1,\mathcal{V}_3)\\
&\le&\sup_{n\in\mathbb{N}}\left\{n^{-2}\mathrm{dist}_H^{D_{L_n}}(CS_n(\mathcal{V}_1),CS_n(\mathcal{V}_3))\right\}\\
&\le&\sup_{n\in\mathbb{N}}\left\{n^{-2}\mathrm{dist}_H^{D_{L_n}}(CS_n(\mathcal{V}_1),CS_n(\mathcal{V}_2))\right.\\
&&+\left.n^{-2}\mathrm{dist}_H^{D_{L_n}}(CS_n(\mathcal{V}_2),CS_n(\mathcal{V}_3))\right\}\\
&\le&\sup_{n\in\mathbb{N}}\left\{n^{-2}\mathrm{dist}_H^{D_{L_n}}(CS_n(\mathcal{V}_1),CS_n(\mathcal{V}_2))\right\}\\
&&+\sup_{n\in\mathbb{N}}\left\{n^{-2}\mathrm{dist}_H^{D_{L_n}}(CS_n(\mathcal{V}_2),CS_n(\mathcal{V}_3))\right\}\\
&\le&\mathrm{dist}_{NC}(CS_n(\mathcal{V}_1),CS_n(\mathcal{V}_2))+\mathrm{dist}_{NC}(CS_n(\mathcal{V}_2),CS_n(\mathcal{V}_3))+2\epsilon.
\end{array}\] Since $\epsilon$ is arbitrary, we obtain
\[\mathrm{dist}_{NC}(\mathcal{V}_1, \mathcal{V}_3)\le
\mathrm{dist}_{NC}(\mathcal{V}_1,
\mathcal{V}_2)+\mathrm{dist}_{NC}(\mathcal{V}_2, \mathcal{V}_3).\]
\end{proof}

\begin{proposition}\label{pro475}
Let $(\mathcal{V},\mathcal{L})$ be a quantized metric space, and let $\mathcal{K}^{(1)}=(K_n^{(1)})$ and 
$\mathcal{K}^{(2)}=(K_n^{(2)})$ 
be compact matrix convex subsets of $\mathcal{CS}(\mathcal{V})$. If $(\mathcal{V}_j,\mathcal{L}_j)$ is the quotient space 
of $(\mathcal{V},\mathcal{L})$ with respect to $\mathcal{K}^{(j)}$ and $\phi^{(j)}$, $j=1,2$, then we have
\[\mathrm{dist}_{NC}(\mathcal{V}_1,\mathcal{V}_2)\le\sup_{k\in\mathbb{N}}\left\{k^{-2}\mathrm{dist}_H^{D_{L_k}}(K^{(1)}_k, 
K^{(2)}_k)\right\}.\]
\end{proposition}

\begin{proof}
Let $p_1$ and $p_2$ be the projections from $\mathcal{V}\oplus\mathcal{V}$ onto the first space $\mathcal{V}$ and the second space 
$\mathcal{V}$, respectively. Denote
\[G^{(j)}_n=(p_1^c)_n(K^{(j)}_n),\ \ \ H^{(j)}_n=(p_2^c)_n(K^{(j)}_n),\ \ \ j=1,2,\ \ n\in\mathbb{N},\]
and set $\mathcal{G}^{(j)}=(G_n^{(j)})$, $\mathcal{H}^{(j)}=(H_n^{(j)})$, $j=1,2$, and $\mathcal{K}=(K_n)=
\overline{\mathrm{mco}}(\mathcal{G}^{(1)}\cup\mathcal{H}^{(2)})$, the BW-closed matrix convex hull of the graded set $(G_n^{(1)}\cup H_n^{(2)})$. Let $(\mathcal{W},\phi(1\oplus 1))$ be the quotient of $(\mathcal{V}\oplus\mathcal{V}, 1\oplus 1)$ with 
respect to $\mathcal{K}$. Then $\mathcal{K}=\phi^c(\mathcal{CS}(\mathcal{W}))$ by Proposition \ref{pro335}. 

For $(a,b)\in\mathrm{Ker}\phi$, we have $\phi(a,b)=0$, that is, $\ll(a,b),\varphi\gg=0_n$ for $\varphi\in K_n$. This is equivalent to 
$\ll(a,b),\varphi\gg=0_n$ for $\varphi\in G_n^{(1)}\cup H_n^{(2)}$, $n\in\mathbb{N}$ since $(a,b)\in A(\mathcal{K})$. And this holds if 
and only if $\ll a,\varphi_1\gg=0_n$ and $\ll b,\varphi_2\gg=0_n$ for $\varphi_1\in G_n^{(1)}$ and $\varphi_2\in H_n^{(2)}$, $n\in\mathbb{N}$, 
that is, if and only if $a\in 
\mathrm{Ker}\phi^{(1)}$ and $b\in\mathrm{Ker}\phi^{(2)}$. So $\mathrm{Ker}\phi=\mathrm{Ker}\phi^{(1)}
\oplus\mathrm{Ker}\phi^{(2)}$. And thus there is a complete order isomorphism $\Psi$ from $\mathcal{W}$ onto $\mathcal{V}_1\oplus\mathcal{V}_2$.

Given $\epsilon>0$. We define a matrix seminorm $\mathcal{N}=(N_n)$ on $\mathcal{V}\oplus\mathcal{V}$ by
\[N_n(a,b)=\epsilon^{-1}\|a-b\|_n,\ \ \ \ a,b\in M_n(\mathcal{V}).\]
Then $\mathcal{N}$ is a matrix bridge between $(\mathcal{V},\mathcal{L})$ and $(\mathcal{V},\mathcal{L})$, and 
$\mathcal{Q}=(Q_n)\in\mathcal{M}(\mathcal{L},\mathcal{L})$ by Proposition \ref{pro52}, where
\[Q_n(a,b)=\max\{L_n(a),L_n(b),N_n(a,b)\},\ \ \ a,b\in M_n(\mathcal{V}), n\in\mathbb{N}.\]
Thus $\mathcal{Q}$ is a matrix Lip-norm on $(\mathcal{V}\oplus\mathcal{V},(1,1))$. Let $\mathcal{P}=(P_n)$ and 
$(\mathcal{W},\mathcal{P})$ be the quotient space of $(\mathcal{V}\oplus\mathcal{V},\mathcal{Q})$ with respect 
to $\mathcal{K}$ and $\phi$. Then $\mathcal{P}\in\mathcal{M}(\mathcal{L}_1,\mathcal{L}_2)$ by Propositon \ref{pro35}. 

Since $D_{P_k}(\varphi_1,\varphi_2)=D_{Q_k}(\phi^c_k(\varphi_1), \phi^c_k(\varphi_2))$ for $\varphi_1, \varphi_2\in CS_k(\mathcal{W})$, 
we have that $\mathrm{dist}_H^{D_{P_k}}(CS_k(\mathcal{V}_1), CS_k(\mathcal{V}_2))=\mathrm{dist}_H^{D_{Q_k}}(G_k^{(1)}, H_k^{(2)})$. For 
$\psi\in K_k^{(2)}$, we have 
\[\begin{array}{rcl}
&&D_{Q_k}((p^c_1)_k(\psi), (p^c_2)_k(\psi))\\
&=&\sup\{\|\ll(p^c_1)_k(\psi), (a,b)\gg-\ll(p^c_2)_k(\psi),(a,b)\gg\|: Q_r(a,b)\le 1, r\in\mathbb{N}\}\\
&\le&\sup\{\|\ll\psi, a\gg-\ll\psi,b\gg\|: N_r(a,b)\le 1, r\in\mathbb{N}\}\\
&\le&\epsilon,\end{array}\]
that is, $\mathrm{dist}_H^{D_{Q_k}}(G_k^{(2)}, H_k^{(2)})\le\epsilon$. Because $\mathcal{Q}\in\mathcal{M}(\mathcal{L},\mathcal{L})$, we get 
that 
\[\mathrm{dist}_H^{D_{Q_k}}(G_k^{(1)}, G_k^{(2)})=\mathrm{dist}_H^{D_{L_k}}(K_k^{(1)}, K_k^{(2)}).\]
So
\[\begin{array}{rcl}
\mathrm{dist}_{NC}(\mathcal{V}_1,\mathcal{V}_2)
&\le&\sup_{k\in\mathbb{N}}\left\{k^{-2}\mathrm{dist}_H^{D_{P_k}}(CS_k(\mathcal{V}_1),CS_k(\mathcal{V}_2))\right\}\\
&=&\sup_{k\in\mathbb{N}}\left\{k^{-2}\mathrm{dist}_H^{D_{Q_k}}(G^{(1)}_k, H^{(2)}_k)\right\}\\
&=&\sup_{k\in\mathbb{N}}\left\{k^{-2}\mathrm{dist}_H^{D_{Q_k}}(G^{(1)}_k, G^{(2)}_k)\right\}\\
&&+\sup_{k\in\mathbb{N}}\left\{k^{-2}\mathrm{dist}_H^{D_{Q_k}}(G^{(2)}_k, H^{(2)}_k)\right\}\\
&\le&\sup_{k\in\mathbb{N}}\left\{k^{-2}\mathrm{dist}_H^{D
Q_{P_k}}(G^{(1)}_k, G^{(2)}_k)+k^{-2}\epsilon\right\}\\
&\le&\sup_{k\in\mathbb{N}}\left\{k^{-2}\mathrm{dist}_H^{D_{Q_k}}(G^{(1)}_k, G^{(2)}_k)\right\}+\epsilon\\
&=&\sup_{k\in\mathbb{N}}\left\{k^{-2}\mathrm{dist}_H^{D_{L_k}}(K^{(1)}_k, K^{(2)}_k)\right\}+\epsilon
\end{array}\]
Since $\epsilon$ is arbitrary, we obtain the desired inequality.
\end{proof}

Let $(A,L_A)$ and $(B,L_B)$ be compact quantum metric spaces. The quantum Gromov-Hausdorff distance between them is defined by
\[\mathrm{dist}_q(A,B)=\inf\mathrm{dist}_H^{\rho_L}(S(A),S(B)),\]
where the infimum is taken over all Lip-norms $L$ on $A\oplus B$ which induce $L_A$ and $L_B$ (see Definition 4.2 in \cite{ri4}).

\begin{proposition}\label{pro49}
Let $(A_j,L_j)$ for $j=1,2$ be compact quantum metric spaces, and let $(\mathcal{V}_j,\mathcal{L}_j)$ be an associated quantized metric space 
of $(A_j,L_j)$(see Example \ref{ex34}). Then
\[\mathrm{dist}_q(A_1, A_2)\le\mathrm{dist}_{NC}(\mathcal{V}_1,\mathcal{V}_2).\]
\end{proposition}

\begin{proof}
Suppose $\mathcal{Q}\in\mathcal{M}(\mathcal{L}_1,\mathcal{L}_2)$. Then $Q_{\mathcal{V}_j,1}=L_{j,1}$for $j=1,2$ and $L_{j,1}(a)=L^s_j(a)$ 
for $a\in A_j$, where $L^s_j=L_{\rho_{L_j}}$ (see Example \ref{ex34}). So for $a\in A_1$, we have
\[\begin{array}{rcl}
L^s_1(a)&=&L_{1,1}(a)=Q_{\mathcal{V}_j,1}(a)\\
&=&\inf\{Q_1(a_1,b_1): \pi_1(a_1,b_1)=a, (a_1,b_1)=\mathcal{V}_1\oplus\mathcal{V}_2\}\\
&=&\inf\{Q_1(a,b_1): b_1\in\mathcal{V}_2\}\\
&\le&\inf\{Q_1(a,b): b\in A_2\}\\
&=&\inf\{R(a,b): b\in A_2\}\\
&=&R_{A_1}(a),
\end{array}\]
where $\pi_1$ is the projection from $\mathcal{V}_1\oplus\mathcal{V}_2$ onto $\mathcal{V}_1$ and $R$ is the restriction of $Q_1$ to $A_1\oplus A_2$. Denote $c=\inf\{Q_1(a,b): b\in\mathcal{V}_2\}$. Let $\epsilon>0$ be given. 
Then there is a $y\in\mathcal{V}_2$ such that $Q_1(a,y)\le c+\epsilon$. Setting $x=\frac12(y+y^\ast)$, we have that $x\in A_2$ and 
\[\begin{array}{rcl}
R(a,x)&=&Q_1(a,x)\\
&=&Q_1\left(a,\frac12(y+y^\ast)\right)\\
&\le&\frac12Q_1((a,y)+(a,y)^\ast)\\
&\le&\frac12(Q_1(a,y)+Q_1((a,y)^\ast))\\
&=&Q_1(a,y)\\
&\le& c+\epsilon.
\end{array}\]
Thus $L^s_1(a)=R_{A_1}(a)$ for $a\in A_1$. Similarly, we have that $L_2^s(b)=R_{A_2}(b)$ for $b\in A_2$. So $R\in\mathcal{M}(L^s_1,L^s_2)$. 

For $\varphi\in CS_1(\mathcal{V}_1)$ and $\psi\in CS_1(\mathcal{V}_2)$, let $\varphi_1=\varphi|_{A_1}$ and $\psi_1=\psi|_{A_2}$. Then 
$\varphi_1\in S(A_1)$ and $\psi_1\in S(A_2)$. Since $Q_1((a,b)^\ast)=Q_1(a,b)$, we obtain
\[\begin{array}{rcl}
D_{Q_1}(\varphi, \psi)&=&\sup\{|\varphi(c)-\psi(d)|: Q_1(c,d)\le 1, (c,d)\in\mathcal{V}_1\oplus\mathcal{V}_2\}\\
&=&\sup\{|\varphi(c)-\psi(d)|: Q_1(c,d)\le 1, (c,d)=(c,d)^\ast\in\mathcal{V}_1\oplus\mathcal{V}_2\}\\
&=&\sup\{|\varphi_1(c)-\psi_1(d)|: R(c,d)\le 1, (c,d)\in A_1\oplus A_2\}\\
&=&\rho_R(\varphi_1,\psi_1),
\end{array}\]
(see \S 2 in \cite{ri4}). So
\[\mathrm{dist}_H^{\rho_R}(S(A_1), S(A_2))=\mathrm{dist}_H^{D_{Q_1}}(CS_1(\mathcal{V}_1),CS_1(\mathcal{V}_2)).\]
Therefore, by Theorem 4.3 and Proposition 7.1 in \cite{ri4}, we have
\[\begin{array}{rcl}
\mathrm{dist}_q(A_1,A_2)&\le&\mathrm{dist}_q((A_1,L_1),(A_1,L_1^s))+\mathrm{dist}_q((A_1,L_1^s),(A_2,L_2^s))\\
&&+\mathrm{dist}_q((A_2,L_2^s),(A_2,L_2))\\
&=&\mathrm{dist}_q((A_1,L_1^s),(A_2,L_2^s))\\
&\le&\mathrm{dist}_H^{\rho_R}(S(A_1),S(A_2))\\
&=&\mathrm{dist}_H^{D_{Q_1}}(CS_1(\mathcal{V}_1), CS_1(\mathcal{V}_2))\\
&\le&\sup_{n\in\mathbb{N}}\left\{n^{-2}\mathrm{dist}_H^{D_{Q_n}}(CS_n(\mathcal{V}_1),CS_n(\mathcal{V}_2))\right\}.
\end{array}\]
Consequently, $\mathrm{dist}_q(A_1,A_2)\le\mathrm{dist}_{NC}(\mathcal{V}_1,\mathcal{V}_2)$.
\end{proof}

\section{Distance zero}\label{s5}

In this section, we show that $\mathrm{dist}_{NC}(\mathcal{V}_1,\mathcal{V}_2)=0$ is equivalent to the existence of a complete 
isometry between them in the following sense.

If $(\mathcal{V}, \mathcal{L})$ is a quantized metric space,
then $\mathcal{L_{D_L}}$ is the largest lower semicontinuous
matrix Lip-norm smaller than $\mathcal{L}$ by Corollary 4.5 in
\cite{wu2}. From Proposition 7.1 in \cite{wu2},
$\mathcal{L_{D_L}}$ extends uniquely to a closed matrix Lip-norm
$\mathcal{L}^c$ on the subspace $\mathcal{V}^c=\{a\in\bar{\mathcal{V}}: L_1^c(a)<+\infty\}$, where $\bar{\mathcal{V}}$  is the 
completion of $\mathcal{V}$ for its matrix norm.

\begin{definition}\label{def53}
Let $(\mathcal{V}_1, \mathcal{L}_1)$ and $(\mathcal{V}_2,
\mathcal{L}_2)$ be quantized metric spaces. By a {\it complete
isometry} from $(\mathcal{V}_1, \mathcal{L}_1)$ onto
$(\mathcal{V}_2, \mathcal{L}_2)$ we mean a unital complete order
isomorphism $\Phi$ from $\mathcal{V}_1^c$ onto $\mathcal{V}_2^c$
such that $\mathcal{L}_1^c=\mathcal{L}_2^c\circ\Phi$, that is,
$L^c_{1,n}=L^c_{2,n}\circ\Phi_n$ for all $n\in\mathbb{N}$.
\end{definition}

\begin{lemma}\label{le54}
Let $(\mathcal{V}, \mathcal{L})$ be a quantized metric space. Then
\[\mathrm{dist}_{NC}(\mathcal{V},\mathcal{V}^c)=0,\ \ \ \ 
\mathrm{dist}_{NC}((\mathcal{V},\mathcal{L}),(\mathcal{V},\mathcal{L_{D_L}}))=0.\]
\end{lemma}

\begin{proof}
Let $\epsilon>0$ be given, and define
\[N_n(a_1,a_2)=\epsilon^{-1}\|a_1-a_2\|_n,\]
for $a_1\in M_n(\mathcal{V})$, $a_2\in M_n(\mathcal{V}^c)$ and
$n\in\mathbb{N}$. Clearly $\mathcal{N}=(N_n)$ is a matrix
continuous matrix seminorm on $\mathcal{V}\oplus\mathcal{V}^c$,
and $N_1(1,1)=0$ and $N_1(1,0)=\epsilon^{-1}\neq 0$.

For $n\in\mathbb{N}$ and $a_1\in M_n(\mathcal{V})$ and
$\delta>0$, setting $a_2=a_1\in M_n(\mathcal{V}^c)$, we have
\[\max\{L_n^c(a_2), N_n(a_1,a_2)\}=L_n^c(a_2)=L_{D_{L_n}}(a_2)\le L_n(a_2)<L_n(a_2)+\delta\]
by Proposition 3.6, Proposition 7.1 and Proposition 3.4 in
\cite{wu2}. Given $n\in\mathbb{N}$ and $a_2\in M_n(\mathcal{V}^c)$
and $\delta>0$. By Lemma 7.3 in \cite{wu2}, there is a sequence
$\{a_1^{(k)}\}$ of elements in $M_n(\mathcal{V})$ such that
$L_n\left(a_1^{(k)}\right)\le L_n^c(a_2)$ and $\{a_1^{(k)}\}$ converges to
$a_2$ in norm. Consequently, we can find an $a_1^{(k_0)}$ such
that $\epsilon^{-1}\left\|a_1^{(k_0)}-a_2\right\|_n\le L_n^c(a_2)+\delta$. So
$\max\left\{L_n\left(a_1^{(k_0)}\right), N_n\left(a_1^{(k_0)},a_2\right)\right\}\le
L_n^c(a_2)+\delta$. Thus $\mathcal{N}$ is a matrix bridge between
$(\mathcal{V},\mathcal{L})$ and $(\mathcal{V}^c,\mathcal{L}^c)$.

Define
\[L_n(a_1,a_2)=\max\{L_n(a_1),L_n^c(a_2),N_n(a_1,a_2)\}\]
for $a_1\in M_n(\mathcal{V})$, $a_2\in M_n(\mathcal{V}^c)$ and
$n\in\mathbb{N}$. By Proposition \ref{pro52},
$\mathcal{L}=(L_n)\in\mathcal{M}(\mathcal{V},\mathcal{V}^c)$. For
$n\in\mathbb{N}$ and $\varphi\in CS_n(\mathcal{V}^c)$, we have
that $\psi=\varphi|_{\mathcal{V}}\in CS_n(\mathcal{V})$, and hence
\[\begin{array}{rcl}
&&D_{L_n}(\psi, \varphi)\\
&=&\sup\big\{\|\ll(\pi_1)_n^c(\psi),
(a_1,a_2)\gg-\ll(\pi_2)_n^c(\varphi), (a_1,a_2)\gg\|:\\
&&L_r(a_1,a_2)\le 1, (a_1,a_2)\in M_r(\mathcal{V}\oplus\mathcal{V}^c), r\in\mathbb{N}\big\}\\
&=&\sup\{\|\ll\psi, a_1\gg-\ll\varphi, a_2\gg\|:\\
&&L_r(a_1,a_2)\le 1, (a_1,a_2)\in M_r(\mathcal{V}\oplus\mathcal{V}^c), r\in\mathbb{N}\}\\
&=&\sup\{\|\ll\varphi, a_1-a_2\gg\|:\\
&&L_r(a_1,a_2)\le 1, (a_1,a_2)\in M_r(\mathcal{V}\oplus\mathcal{V}^c), r\in\mathbb{N}\}\\
&\le&\sup\{\|a_1-a_2\|_r:
L_r(a_1,a_2)\le 1, (a_1,a_2)\in M_r(\mathcal{V}\oplus\mathcal{V}^c), r\in\mathbb{N}\}\\
&\le&\epsilon,
\end{array}\] 
where $\pi_1$ and $\pi_2$ are the projections from $\mathcal{V}\oplus\mathcal{V}^c$ onto $\mathcal{V}$ and $\mathcal{V}^c$, respectively. 
For $n\in\mathbb{N}$ and $\varphi\in CS_n(\mathcal{V})$,
there is a $\psi\in CS_n(\mathcal{V}^c)$ such that
$\psi|_{\mathcal{V}}=\varphi$ by Arveson's extension theorem. So
\[\begin{array}{rcl}
&&D_{L_n}(\varphi, \psi)\\
&=&\sup\big\{\|\ll(\pi_1)_n^c(\varphi),
(a_1,a_2)\gg-\ll(\pi_2)_n^c(\psi), (a_1,a_2)\gg\|:\\
&&L_r(a_1,a_2)\le 1, (a_1,a_2)\in M_r(\mathcal{V}\oplus\mathcal{V}^c), r\in\mathbb{N}\big\}\\
&=&\sup\{\|\ll\varphi, a_1\gg-\ll\psi, a_2\gg\|:\\
&&L_r(a_1,a_2)\le 1, (a_1,a_2)\in M_r(\mathcal{V}\oplus\mathcal{V}^c), r\in\mathbb{N}\}\\
&=&\sup\{\|\ll\psi, a_1-a_2\gg\|:\\
&&L_r(a_1,a_2)\le 1, (a_1,a_2)\in M_r(\mathcal{V}\oplus\mathcal{V}^c), r\in\mathbb{N}\}\\
&\le&\sup\{\|a_1-a_2\|_r:
L_r(a_1,a_2)\le 1, (a_1,a_2)\in M_r(\mathcal{V}\oplus\mathcal{V}^c), r\in\mathbb{N}\}\\
&\le&\epsilon.
\end{array}\]
Thus
$\mathrm{dist}_H^{D_{L_n}}(CS_n(\mathcal{V}),CS_n(\mathcal{V}^c))\le\epsilon$
for $n\in\mathbb{N}$, and so
\[\sup_{n\in\mathbb{N}}\left\{\mathrm{dist}_H^{D_{L_n}}(CS_n(\mathcal{V}),CS_n(\mathcal{V}^c))\right\}\le\epsilon.\]
Therefore,
$\mathrm{dist}_{NC}(\mathcal{V},\mathcal{V}^c)\le\epsilon$. By the
arbitrariness of $\epsilon$, we obtain
\[\mathrm{dist}_{NC}(\mathcal{V},\mathcal{V}^c)=0.\]

By Proposition 3.4 in \cite{wu2} and the proof of Theorem 4.4 in \cite{wu2}, we can prove that 
$\mathrm{dist}_{NC}((\mathcal{V},\mathcal{L}),(\mathcal{V},\mathcal{L_{D_L}}))=0$ similarly.
\end{proof}

\begin{theorem}\label{th55}
Suppose $(\mathcal{V}_1, \mathcal{L}_1)$ and $(\mathcal{V}_2,
\mathcal{L}_2)$ are quantized metric spaces. If there exists
a complete isometry $\Phi$ from $(\mathcal{V}_1, \mathcal{L}_1)$
onto $(\mathcal{V}_2, \mathcal{L}_2)$, then
\[\mathrm{dist}_{NC}(\mathcal{V}_1, \mathcal{V}_2)=0.\]
\end{theorem}

\begin{proof}
For $\epsilon>0$, we define
\[N_n(a_1,a_2)=\epsilon^{-1}\|\Phi_n(a_1)-a_2\|_n,\]
for $a_1\in M_n(\mathcal{V}_1^c)$, $a_2\in M_n(\mathcal{V}_2^c)$
and $n\in\mathbb{N}$. Clearly $\mathcal{N}=(N_n)$ is a matrix
seminorm on $\mathcal{V}_1^c\oplus\mathcal{V}_2^c$ and $N_n((a_1,a_2)^\ast)=N_n(a_1,a_2)$
for $a_1\in M_n(\mathcal{V}_1^c)$, $a_2\in M_n(\mathcal{V}_2^c)$
and $n\in\mathbb{N}$. And we have
that $N_1(1,1)=\epsilon^{-1}\|\Phi(1)-1\|_1=0$ and  
$N_1(1,0)=\epsilon^{-1}\|\Phi(1)-0\|_1=\epsilon^{-1}$. If
$\{a_1^{(k)}\}\subseteq M_n(\mathcal{V}_1^c)$ and
$\{a_2^{(k)}\}\subseteq M_n(\mathcal{V}_2^c)$ with
$\lim_{k\to\infty}a_1^{(k)}=a_1\in M_n(\mathcal{V}_1^c)$ and
$\lim_{k\to\infty}a_2^{(k)}=a_2\in M_n(\mathcal{V}_2^c)$, we have
that $\lim_{k\to\infty}
N_n(a_1^{(k)},a_2^{(k)})=\lim_{k\to\infty}\epsilon^{-1}\|\Phi_n(a_1^{(k)})-a_2^{(k)}\|_n=\epsilon^{-1}
\|\Phi_n(a_1)-a_2\|_n=N_n(a_1,a_2)$ since $\Phi$ is completely
bounded(see Proposition 3.5 in \cite{pau}).

Given $a_1\in M_n(\mathcal{V}_1^c)$ and $\delta>0$. Taking
$a_2=\Phi_n(a_1)$, we have
\[\begin{array}{rcl}\max\{L^c_{2,n}(a_2),N_n(a_1,a_2)\}&=&\max\{L^c_{2,n}(\Phi_n(a_1)),
\epsilon^{-1}\|\Phi_n(a_1)-a_2\|\}\\
&=&L^c_{1,n}(a_1)<L^c_{1,n}(a_1)+\delta.\end{array}\] While if
$a_2\in M_n(\mathcal{V}_2^c)$ and $\delta>0$, we can take $a_1\in
M_n(\mathcal{V}_1^c)$ such that $\Phi_n(a_1)=a_2$, and hence we
have
\[\begin{array}{rcl}\max\{L^c_{1,n}(a_1),N_n(a_1,a_2)\}&=&\max\{L^c_{2,n}(\Phi_n(a_1)),
\epsilon^{-1}\|\Phi_n(a_1)-a_2\|\}\\
&=&L^c_{2,n}(a_2)<L^c_{2,n}(a_2)+\delta.\end{array}\] Therefore,
$\mathcal{N}$ is a matrix bridge between $(\mathcal{V}_1^c,
\mathcal{L}_1^c)$ and $(\mathcal{V}_2^c, \mathcal{L}_2^c)$.

Define
\[L_n(a_1,a_2)=\max\{L^c_{1,n}(a_1), L^c_{2,n}(a_2), N_n(a_1,a_2)\},\]
for $a_1\in M_n(\mathcal{V}_1^c)$, $a_2\in M_n(\mathcal{V}_2^c)$
and $n\in\mathbb{N}$. By Proposition \ref{pro52},
$\mathcal{L}=(L_n)\in\mathcal{M}(\mathcal{L}_1^c,\mathcal{L}_2^c)$.
For $n\in\mathbb{N}$ and $\varphi\in CS_n(\mathcal{V}_2^c)$, we
have that $\varphi\circ\Phi\in CS_n(\mathcal{V}_1^c)$, and so
\[\begin{array}{rcl}
&&D_{L_n}(\varphi\circ\Phi, \varphi)\\
&=&\sup\{\|\ll\varphi\circ\Phi, a_1\gg-\ll\varphi, a_2\gg\|:\\
&&L_r(a_1,a_2)\le 1, (a_1,a_2)\in M_r(\mathcal{V}_1^c\oplus\mathcal{V}_2^c), r\in\mathbb{N}\}\\
&=&\sup\{\|\ll\varphi, \Phi_r(a_1)-a_2\gg\|:\\
&&L_r(a_1,a_2)\le 1, (a_1,a_2)\in M_r(\mathcal{V}_1^c\oplus\mathcal{V}_2^c), r\in\mathbb{N}\}\\
&\le&\sup\{\|\Phi_r(a_1)-a_2\|_r:
L_r(a_1,a_2)\le 1, (a_1,a_2)\in M_r(\mathcal{V}_1^c\oplus\mathcal{V}_2^c), r\in\mathbb{N}\}\\
&\le&\epsilon,
\end{array}\]
Similarly, for $n\in\mathbb{N}$ and $\psi\in
CS_n(\mathcal{V}_1^c)$, we have that $D_{L_n}((\psi,\psi\circ\Phi^{-1})\le\epsilon$. Thus we obtain that
$\mathrm{dist}_H^{D_{L_n}}(CS_n(\mathcal{V}_1^c),CS_n(\mathcal{V}_2^c))\le\epsilon$
for $n\in\mathbb{N}$, and so
\[\sup_{n\in\mathbb{N}}\left\{n^{-2}\mathrm{dist}_H^{D_{L_n}}(CS_n(\mathcal{V}_1^c),CS_n(\mathcal{V}_2^c))\right\}\le\epsilon.\]
Therefore, $\mathrm{dist}_{NC}(\mathcal{V}_1^c,
\mathcal{V}_2^c)\le\epsilon$. Since $\epsilon$ is arbitrary, we
conclude 
$$\mathrm{dist}_{NC}(\mathcal{V}_1^c,\mathcal{V}_2^c)=0.$$
Now, by Theorem \ref{th42} and Lemma \ref{le54} we have
\[\mathrm{dist}_{NC}(\mathcal{V}_1,\mathcal{V}_2)\le\mathrm{dist}_{NC}(\mathcal{V}_1,\mathcal{V}_1^c)
+\mathrm{dist}_{NC}(\mathcal{V}_1^c,\mathcal{V}_2^c)+\mathrm{dist}_{NC}(\mathcal{V}_2^c,\mathcal{V}_2)=0.\]
So $\mathrm{dist}_{NC}(\mathcal{V}_1,\mathcal{V}_2)=0$.\end{proof}

Given a quantized metric space $(\mathcal{V}, \mathcal{L})$. From Proposition 6.1 in \cite{wu2} and the proof of Proposition 3.5 in 
\cite{wewi}, the mapping $\Psi: \mathcal{V}\mapsto A(\mathcal{CS}(\mathcal{V}))$, defined by $\Psi(a)(\varphi)=\varphi(a)$ for 
$\varphi\in CS_r(\mathcal{V})$, is a unital complete order isomorphism from $\mathcal{V}$ into
$A(\mathcal{CS(V)})$, and $\Psi$ can be extended to a unital complete order isomorphism $\bar{\Psi}$ from the completion $\bar{\mathcal{V}}$ of $\mathcal{V}$ onto
$A(\mathcal{CS(V)})$. Define
\[L_{D_{L_n}}(\mathbf{F}^{(n)})=\sup\left\{\frac{\| F_r^{(n)}(\varphi)-F_r^{(n)}(\psi)\|} {D_{L_r}(\varphi,
\psi)}: \varphi\neq\psi, \varphi, \psi\in CS_r(\mathcal V), r\in\mathbb{N}\right\},\] 
where $\mathbf F^{(n)}\in A(\mathcal{CS(V)},
M_n)$. Then $\mathcal{L_{D_L}}=(L_{D_{L_n}})$ is a matrix gauge on
$A(\mathcal{CS(V)})$. Denote
\[K_n=\left\{\mathbf{F}^{(n)}\in A(\mathcal{CS(V)}, M_n): L_{D_{L_n}}(\mathbf{F}^{(n)})<+\infty
\right\}.\]
Let $L_n^1=\{a\in\mathcal{V}:L_n(a)\le 1\}$ and $\bar{L}_n^1$ be the norm closure of $L_n^1$ in $\bar{\mathcal{V}}$. Denote $\mathcal{L}^1=
(L_n^1)$ and $\bar{\mathcal{L}}^1=(\bar{L}_n^1)$. The matrix gauge $\bar{\mathcal{L}}=(\bar{L}_n)$ on $(\bar{\mathcal{V}},1)$ determined by 
$\bar{\mathcal{L}}^1$ is called the {\it closure} of $\mathcal{L}$. $\mathcal{L}$ is {\it closed} if $\mathcal{L}=\bar{\mathcal{L}}$ on the subspace 
where $\bar{\mathcal{L}}$ is finite (see Definition 7.2 in \cite{wu2}).

\begin{lemma}\label{le56}
If $\mathcal{L}$ is closed, then $\Psi_n(M_n(\mathcal{V}))=K_n$
for $n\in\mathbb{N}$.
\end{lemma}

\begin{proof}
Denote
\[M_n^1=\{\Psi_n(a): a\in M_n(\mathcal{V}), L_n(a)\le 1\}, n\in\mathbb{N},\]
\[L_{D_{L_n}}^1=\left\{\mathbf{F}^{(n)}\in A(\mathcal{CS(V)}, M_n): L_{D_{L_n}}(\mathbf{F}^{(n)})\le 1\right\}, n\in\mathbb{N},\]
and set $\mathcal{M}^1=(M_n^1)$. Define
\[L_n^{\prime}(f)=\sup\left\{\|\ll f, \Psi_r(a)\gg\|: L_r(a)\le 1, a\in M_r(\mathcal{V}), r\in\mathbb{N}\right\},\]
for $f\in M_n((A(\mathcal{CS}(\mathcal{V}))/(\mathbb{C}\mathbf{I}))^*)$, where $\mathbf{I}$ is the order unit of 
$A(\mathcal{CS}(\mathcal{V}))$. Here we view 
$M_n((A(\mathcal{CS}(\mathcal{V}))/(\mathbb{C}\mathbf{I}))^\ast)$ as the subspace of 
$M_n((A(\mathcal{CS}(\mathcal{V})))^*)$ consisting of those 
$f\in M_n((A(\mathcal{CS}(\mathcal{V})))^\ast)$ with $f(a)=0_n$ for $a\in\mathbb{C}\mathbf{I}$. 
Clearly, $\mathcal{L}^\prime=(L_n^\prime)$ is a matrix gauge on $(A(\mathcal{CS}(\mathcal{V}))/(\mathbb{C}\mathbf{I}))^\ast$ and 
$L_n^\prime(f^\ast)=L_n^\prime(f)$ for all $f\in M_n((A(\mathcal{CS}(\mathcal{V}))/(\mathbb{C}\mathbf{I}))^*)$,
The generalized bipolar theorem says that $({\mathcal M}^1)^{\circledcirc\circledcirc}$ is
the smallest weakly closed absolutely matrix convex set containing
$\mathcal M^1$ (see Proposition 4.1 in \cite{efwe}). Since ${\mathcal L}=(L_n)$ is a matrix
gauge and $\Psi$ is a unital complete order isomorphism,
${\mathcal M}^1$ is absolutely matrix convex. The closeness of
$\mathcal L$ implies that $\mathcal M^1$ is normed-closed by Lemma 7.4 in \cite{wu2}. So $\mathcal{M}^1$ is weakly closed. Thus
\[({\mathcal M}^1)^{\circledcirc\circledcirc}=\mathcal M^1.\]
For $n\in\mathbb N$, we have
\[\begin{array}{rcl}
(M_n^1)^{\circledcirc}
&=&\{\Psi_n(a): a\in M_n(\mathcal V), L_n(a)\le 1\}^\circledcirc\\
&=&\{f\in M_n((A(\mathcal{CS(V)}))^*):\,\|\ll f, \Psi_r(a)\gg\|\le
1\\
&&\hbox{ for all }a\in M_r(\mathcal V),
L_r(a)\le 1, r\in{\mathbb N}\}\\
&=&\{f\in M_n((A(\mathcal{CS}(\mathcal{V}))/(\mathbb{C}\mathbf{I}))^*):\,\|\ll f, \Psi_r(a)\gg\|\le
1\\
&&\hbox{ for all }a\in M_r(\mathcal V), L_r(a)\le 1, r\in{\mathbb
N}\}\\
&=&\{f\in M_n((A(\mathcal{CS}(\mathcal{V}))/(\mathbb{C}\mathbf{I}))^*):\,L_n^{\prime}(f)\le 1\}
\end{array}\]
and
\[\begin{array}{rcl}
&&(M_n^1)^{\circledcirc\circledcirc}\\
&=&\{f\in M_n((A(\mathcal{CS}(\mathcal{V}))/(\mathbb{C}\mathbf{I}))^*):\,L_n^{\prime}(f)\le 1\}^\circledcirc\\
&=&\{\mathbf{F}^{(n)}\in M_n(A(\mathcal{CS(V)})):\,\|\ll f,\mathbf{F}^{(n)}\gg\|\le 1\\
&&\hbox{ for all }f\in M_r((A(\mathcal{CS}(\mathcal{V}))/(\mathbb{C}\mathbf{I}))^*), L^{\prime}_r(f)\le
1, r\in{\mathbb N}\}.
\end{array}\]
So $\mathbf{F}^{(n)}\in (M^1_n)^{\circledcirc\circledcirc}$ if and only if 
$$\|\ll f, \mathbf{F}^{(n)}\gg\|\le L_r^\prime(f)$$
for all $f\in M_r((A(\mathcal{CS}(\mathcal{V}))/(\mathbb{C}\mathbf{I}))^\ast)$ and $r\in\mathbb{N}$. 

Suppose that 
$\|\ll f, \mathbf{F}^{(n)}\gg\|\le L_r^\prime(f)$ for all $f=f^\ast\in M_r((A(\mathcal{CS}(\mathcal{V}))/(\mathbb{C}\mathbf{I}))^\ast)$ 
and $r\in\mathbb{N}$. Then for $g\in M_r((A(\mathcal{CS}(\mathcal{V}))/(\mathbb{C}\mathbf{I}))^\ast)$, we have
\[\begin{array}{rcl}
\|\ll g, \mathbf{F}^{(n)}\gg\|
&=&\left\|\left[\begin{array}{cc}
1&0
\end{array}\right]
\left[\begin{array}{cc}
0&\ll g,\mathbf{F}^{(n)}\gg\\
\ll g^\ast,\mathbf{F}^{(n)}\gg&0
\end{array}\right]\left[\begin{array}{cc}
0\\
1\end{array}\right]\right\|\\
&\le&\left\|\ll\left[\begin{array}{cc}
0&g\\
g^\ast&0\end{array}\right], \mathbf{F}^{(n)}\gg\right\|\\
&\le&L^\prime_{2r}\left(\left[\begin{array}{cc}
0&g\\
g^\ast&0\end{array}\right]\right)\\
&\le&L^\prime_{2r}\left(\left[\begin{array}{cc}
g&0\\
0&g^\ast\end{array}\right]\left[\begin{array}{cc}
0&1\\
1&0\end{array}\right]\right)\\
&\le&L^\prime_{2r}\left(\left[\begin{array}{cc}
g&0\\
0&g^\ast\end{array}\right]\right)\\
&=&L_r^\prime(g),\end{array}\] 
Thus $\mathbf{F}^{(n)}\in (M^1_n)^{\circledcirc\circledcirc}$ exactly if 
$\|\ll f, \mathbf{F}^{(n)}\gg\|\le L_r^\prime(f)$ for all $f=f^\ast\in M_r((A(\mathcal{CS}(\mathcal{V}))/(\mathbb{C}\mathbf{I}))^\ast)$ 
and $r\in\mathbb{N}$. According to Lemma 4.1 in \cite{wu2}, $\mathbf{F}^{(n)}\in (M^1_n)^{\circledcirc\circledcirc}$ exactly if 
$\|\ll\varphi, \mathbf{F}^{(n)}\gg-\ll\psi, \mathbf{F}^{(n)}\gg\|\le L_r^\prime(\varphi-\psi)$ for all 
$\varphi,\psi\in CS_r(A(\mathcal{CS}(\mathcal{V})))$ and $r\in\mathbb{N}$. So $\mathbf{F}^{(n)}\in (M^1_n)^{\circledcirc\circledcirc}$ 
exactly if 
\[\begin{array}{rcl}
&&\|\ll\varphi\circ\Psi^{-1}, \mathbf{F}^{(n)}\gg-\ll\psi\circ\Psi^{-1}, \mathbf{F}^{(n)}\gg\|\\
&\le&L_r^\prime(\varphi\circ\Psi^{-1}-\psi\circ\Psi^{-1})\\
&=&\sup\left\{\|\ll(\varphi-\psi)\circ\Psi^{-1}, \Psi_k(a)\gg\|:L_k(a)\le 1, a\in M_k(\mathcal{V}), k\in\mathbb{N}\right\}\\
&=&\sup\{\|\ll\varphi-\psi, a\gg\|:L_k(a)\le 1, a\in M_k(\mathcal{V}), k\in\mathbb{N}\}\\
&=&D_{L_r}(\varphi,\psi).\end{array}\]
for all $\varphi,\psi\in CS_r(A(\mathcal{CS}(\mathcal{V})))$ and $r\in\mathbb{N}$. Because 
$\|\ll\varphi\circ\Psi^{-1}, \mathbf{F}^{(n)}\gg-\ll\psi\circ\Psi^{-1}, \mathbf{F}^{(n)}\gg\|=\|F_r^{(n)}(\varphi)-F_r^{(n)}(\psi)\|$, 
$\mathbf{F}^{(n)}\in (M^1_n)^{\circledcirc\circledcirc}$ if and only if $\|F_r^{(n)}(\varphi)-F_r^{(n)}(\psi)\|\le D_{L_r}(\varphi,\psi)$ 
for all $\varphi,\psi\in CS_r(A(\mathcal{CS}(\mathcal{V})))$ and $r\in\mathbb{N}$. And this says exactly that 
$\mathbf{F}^{(n)}\in L^1_{D_{L_n}}$. Therefore, $M_n^1=L^1_{D_{L_n}}$ for $n\in\mathbb{N}$. So 
$\Psi_n(M_n(\mathcal{V}))=K_n$ for $n\in\mathbb{N}$.
\end{proof}

\begin{lemma}\label{le57}
Let $(\mathcal{V}_1, \mathcal{L}_1)$ and $(\mathcal{V}_2,
\mathcal{L}_2)$ be quantized metric spaces such that
$\mathcal{L}_1$ and $\mathcal{L}_2$ are closed. For every matrix
affine mapping $\alpha=(\alpha_n)$ from
$\mathcal{CS}(\mathcal{V}_1)$ onto $\mathcal{CS}(\mathcal{V}_2)$
which is completely isometric for $\mathcal{D}_{\mathcal{L}_1}$
and $\mathcal{D}_{\mathcal{L}_2}$, there is a unital complete
order isomorphism $\Psi$ from $\mathcal{V}_1$ onto $\mathcal{V}_2$
such that $\mathcal{L}_1=\mathcal{L}_2\circ\Psi$.
\end{lemma}

\begin{proof}
Define $\Phi: A(\mathcal{CS}(\mathcal{V}_{2}))\longmapsto
A(\mathcal{CS}(\mathcal{V}_{1}))$ by
\[\left(\Phi(\mathbf{F}^{(1)})\right)_r(\varphi)=F^{(1)}_r(\alpha_r(\varphi))\]
for $\mathbf{F}^{(1)}\in A(\mathcal{CS}(\mathcal{V}_{2}))$ and
$\varphi\in CS_r(\mathcal{V}_1)$. Since $\alpha$ is isometric and matrix affine, $\Phi$ is well-defined. Clearly, $\Phi$ is unital and
surjective. On the level of matrices, we have
\[\left(\Phi_n(\mathbf{F}^{(n)})\right)_r(\varphi)=F^{(n)}_r(\alpha_r(\varphi))\]
for $\mathbf{F}^{(n)}\in M_n(A(\mathcal{CS}(\mathcal{V}_{2})))$
and $\varphi\in CS_r(\mathcal{V}_1)$. Since $\mathbf{F}^{(n)}\ge
0$ in $M_n(A(\mathcal{CS}(\mathcal{V}_{2})))$ if and only if
$F^{(n)}_r(\varphi)\ge 0$ for all $r\in\mathbb{N}$ and $\varphi\in
CS_r(\mathcal{V}_2)$, $\Phi$ is a unital complete order
isomorphism from $A(\mathcal{CS}(\mathcal{V}_{2}))$ onto
$A(\mathcal{CS}(\mathcal{V}_{1}))$. Since $\mathcal{L}_1$ and
$\mathcal{L}_2$ are closed, $\Phi$ is a unital complete order
isomorphism from $\mathcal{V}_2$ onto $\mathcal{V}_1$ by Lemma
\ref{le56}. That $\alpha$ is completely isometric for
$\mathcal{D}_{\mathcal{L}_1}$ and $\mathcal{D}_{\mathcal{L}_2}$
implies that $L_{D_{L_{1,n}}}(\Phi_n(a_2))=L_{D_{L_{2,n}}}(a_2)$
for all $a_2\in M_n(\mathcal{V}_2)$ and $n\in\mathbb{N}$. Because
$\mathcal{L}_1$ and $\mathcal{L}_2$ are closed, they are lower
semicontinuous, so that $\mathcal{L}_{\mathcal{D}_{\mathcal{L}_{1}}}=\mathcal{L}_1$ on
$\mathcal{V}_1$ by Theorem 4.4 in \cite{wu2}, and similarly for $\mathcal{V}_2$. Thus
$\Phi^{-1}$ is a unital complete order isomorphism from
$\mathcal{V}_1$ onto $\mathcal{V}_2$ such that
$\mathcal{L}_1=\mathcal{L}_2\circ\Psi$.
\end{proof}

\begin{theorem}\label{th58}
Suppose $(\mathcal{V}_1, \mathcal{L}_1)$ and $(\mathcal{V}_2,
\mathcal{L}_2)$ are quantized metric spaces. If
\[\mathrm{dist}_{NC}(\mathcal{V}_1, \mathcal{V}_2)=0,\]
then there
exists a complete isometry $\Phi$ from $(\mathcal{V}_1,
\mathcal{L}_1)$ onto $(\mathcal{V}_2, \mathcal{L}_2)$.
\end{theorem}

\begin{proof}
Since $\mathrm{dist}_{NC}(\mathcal{V}_1, \mathcal{V}_2)=0$, we
have
\[\mathrm{dist}_{NC}(\mathcal{V}_1^c, \mathcal{V}_2^c)=0\]
by Lemma \ref{le54} and Theorem \ref{th42}. From that
$\mathrm{dist}_{NC}(\mathcal{V}_1^c, \mathcal{V}_2^c)=0$, there is
a sequence
$\{\mathcal{L}^{(k)}\}\subseteq\mathcal{M}(\mathcal{V}_1^c\oplus\mathcal{V}_2^c)$
of matrix Lip-norms such that
\[\sup_{n\in\mathbb{N}}\left\{n^{-2}\mathrm{dist}_H^{D_{L_n^{(k)}}}(CS_n(\mathcal{V}_1^c),
CS_n(\mathcal{V}^c_2))\right\}<\frac{1}{k}.\] So for each
$n\in\mathbb{N}$, we have
\[n^{-2}\mathrm{dist}_H^{D_{L_n^{(k)}}}(CS_n(\mathcal{V}^c_1),
CS_n(\mathcal{V}^c_2))<\frac{1}{k}.\] And for $\varphi,\psi\in
CS_n(\mathcal{V}^c_i)$, $i=1,2$, by Proposition \ref{pro33} we
have
\[D_{L_n^{(k)}}((\pi_i)^c_n(\varphi),(\pi_i)^c_n(\psi))=D_{L^c_{i,n}}(\varphi,\psi),\]
where $\pi_i, i=1,2$, is the projection from $\mathcal{V}_1^c\oplus\mathcal{V}_2^c$ onto $\mathcal{V}_i^c$. 
Therefore, for each $n\in\mathbb{N}$ we get
\[\mathrm{dist}_{GH}(CS_n(\mathcal{V}^c_1),
CS_n(\mathcal{V}^c_2))=0,\] where $\mathrm{dist}_{GH} (CS_n
(\mathcal{V}^c_1), CS_n(\mathcal{V}^c_2))$ is the Gromov-Hausdorff
distance(see Definition 3.4 in \cite{gro}) between $(CS_n(\mathcal{V}^c_1),
D_{L^c_{1,n}})$ and $(CS_n(\mathcal{V}^c_2), D_{L^c_{2,n}})$. As
in the proofs of Theorem 7.6 and Theorem 7.7 in \cite{ri4}, there
is a subsequence $\left\{D_{L_1^{(k_{j_1})}}\right\}$ which
converges uniformly on the disjoint union
$CS_1(\mathcal{V}^c_1)\sqcup CS_1(\mathcal{V}^c_2)$ to a
semi-metric $\sigma_1$ and $\sigma_1$ determines an isometry
$\alpha_1$ from $CS_1(\mathcal{V}^c_1)$ onto
$CS_1(\mathcal{V}^c_2)$ by the condition that
$\sigma_1(\varphi,\alpha_1(\varphi))=0$. Similarly, there is a
subsequence $\left\{D_{L_2^{(k_{j_1,j_2})}}\right\}$ of
$\left\{D_{L_2^{(k_{j_1})}}\right\}$ which converges uniformly on
$CS_2(\mathcal{V}^c_1)\sqcup CS_2(\mathcal{V}^c_2)$ to a
semi-metric $\sigma_2$ and $\sigma_2$ determines an isometry
$\alpha_2$ from $CS_2(\mathcal{V}^c_1)$ onto
$CS_2(\mathcal{V}^c_2)$ by the condition that
$\sigma_2(\varphi,\alpha_2(\varphi))=0$. In general, once
\[\left\{D_{L_2^{(k_{j_1})}}\right\},
\left\{D_{L_2^{(k_{j_1,j_2})}}\right\},\cdots,\left\{D_{L_2^{(k_{j_1,j_2,
\cdots,j_n})}}\right\}\]
 have been choosen, there is a subsequence
$\left\{D_{L_2^{(k_{j_1,j_2,\cdots,j_n,j_{n+1}})}}\right\}$ of
$\left\{D_{L_2^{(k_{j_1,j_2,\cdots,j_n})}}\right\}$ which
converges uniformly on $CS_{n+1}(\mathcal{V}^c_1)\sqcup
CS_{n+1}(\mathcal{V}^c_2)$ to a semi-metric $\sigma_{n+1}$ and
$\sigma_{n+1}$ determines an isometry $\alpha_{n+1}$ from
$CS_{n+1}(\mathcal{V}^c_1)$ onto $CS_{n+1}(\mathcal{V}^c_2)$ by
the condition that
$\sigma_{n+1}(\varphi,\alpha_{n+1}(\varphi))=0$.

Given $\varphi_i\in CS_{n_i}(\mathcal{V}^c_1)$ and $\gamma_i\in
M_{n_i,n}$, $i=1,2,\cdots,m$, satisfying
$\sum_{i=1}^m\gamma_i^*\gamma_i=1_n$, and $\epsilon>0$. Let $s=\max\{n,n_1,\cdots,n_m\}$. We can
find $K\in\mathbb{N}$ such that if $k_{j_1,\cdots,j_s}>K$ then
$\left\|\sigma_l-D_{L_l^{(k_{j_1,\cdots,j_s})}}\right\|<\frac{\epsilon}{2}$ for
$l=n,n_1,\cdots,n_m$ . Now for $k_{j_1,\cdots,j_s}>K$ we have
\[\begin{array}{rcl}
&&\sigma_n\left(\sum_{i=1}^m\gamma_i^*\varphi_i\gamma_i,\sum_{i=1}^m\gamma_i^*\alpha_{n_i}(\varphi_i)\gamma_i\right)\\
&\le&D_{L_n^{(k_{j_1,\cdots,j_s})}}\left(\sum_{i=1}^m\gamma_i^*\varphi_i\gamma_i,
\sum_{i=1}^m\gamma_i^*\alpha_{n_i}(\varphi_i)\gamma_i\right)+\frac{\epsilon}{2}\\
&\le&D_{L_{n_1+\cdots+n_m}^{(k_{j_1,\cdots,j_s})}}\left(\varphi_1\oplus\cdots\oplus\varphi_m,
\alpha_{n_1}(\varphi_1)\oplus\cdots\oplus\alpha_{n_m}(\varphi_m)\right)+\frac{\epsilon}2\\
&=&\max\left\{D_{L_{n_1}^{(k_{j_1,\cdots,j_s})}}(\varphi_1, \alpha_{n_1}(\varphi_1)),\cdots,D_{L_{n_m^{(k_{j_1,\cdots,j_s})}}}(\varphi_m,
\alpha_{n_m}(\varphi_m))\right\}+\frac{\epsilon}2\\
&<&\max\left\{\sigma_{n_1}(\varphi_1, \alpha_{n_1}(\varphi_1))+\frac{\epsilon}2,\cdots,\sigma_{n_m}(\varphi_m,
\alpha_{n_m}(\varphi_m))+\frac{\epsilon}2\right\}+\frac{\epsilon}2\\
&=&\epsilon.
\end{array}\]
Since $\epsilon$ is arbitrary, we have $\sigma_n\left(\sum_{i=1}^m\gamma_i^*\varphi_i\gamma_i,
\sum_{i=1}^m\gamma_i^*\alpha_{n_i}(\varphi_i)\gamma_i\right)=0$. But 
$$\sigma_n\left(\sum_{i=1}^m\gamma_i^*\varphi_i\gamma_i,
\alpha_n\left(\sum_{i=1}^m\gamma_i^*\varphi_i\gamma_i\right)\right)=0.$$
By Lemma 7.4 in \cite{ri4}, we obtain 
\[\alpha_n\left(\sum_{i=1}^m\gamma_i^*\varphi_i\gamma_i\right)=\sum_{i=1}^m\gamma_i^*\alpha_{n_i}(\varphi_i)\gamma_i.\]
So $\alpha=(\alpha_n)$ is matrix affine.

Now, by Lemma \ref{le57} we conclude that there exists a unital
complete order isomorphism $\Phi$ from $\mathcal{V}^c_1$ onto
$\mathcal{V}^c_2$ such that
$\mathcal{L}^c_1=\mathcal{L}_2^c\circ\Phi$, that is, $\Phi$ is a
complete isometry from $(\mathcal{V}_1, \mathcal{L}_1)$ onto
$(\mathcal{V}_2, \mathcal{L}_2)$.
\end{proof}

\section{Completeness}\label{s6}

For the metric space of complete isometry classes of quantized metric spaces with the quantized Gromov-Hausdorff distance, we show in this 
section that it is complete.

Let $\{(\mathcal{V}_i,1)\}$ be a sequence of matrix order unit space. We will denote by $\oplus_{i\in\mathbb{N}}\mathcal{V}_i$ the
operator space direct sum that is formed of all sequences $\{a_i\}$ with $a_i\in\mathcal{V}_i$ and $\sup_{i\in\mathbb{N}}\|a_i\|<
+\infty$, and by $\oplus_{i=1}^n\mathcal{V}_i$ the operator space direct sum of $\mathcal{V}_1,\mathcal{V}_2,\cdots,\mathcal{V}_n$ (see \S2.6 in 
\cite{pi1}). Then $(\oplus_{i\in\mathbb{N}}\mathcal{V}_i, \{1\})$ and $(\oplus_{i=1}^n\mathcal{V}_i,(\underbrace{1,\cdots,1}_n))$ are 
matrix order unit spaces.

Suppose we have a sequence $\{(\mathcal{V}_i,\mathcal{L}_i)\}$ of quantized metric spaces. Suppose further that we have a 
sequence $\{\mathcal{M}_i\}$ of matrix Lip-norms with $\mathcal{M}_i\in\mathcal{M}(\mathcal{L}_i,\mathcal{L}_{i+1})$. Define $\mathcal{Q}=
(Q_k)$ on $\prod_{i\in\mathbb{N}}\mathcal{V}_i$, the full product, by
\[Q_k(\{a_i\})=\sup_{i\in\mathbb{N}}\{M_{i,k}(a_i,a_{i+1})\},\ \ \ \ \{a_i\}\in M_k\left(\prod_{i\in\mathbb{N}}\mathcal{V}_i\right),\]
and set
\[\mathcal{E}_1=\{\{a_i\}\in\oplus_{i\in\mathbb{N}}\mathcal{V}_i: Q_1(\{a_i\})<+\infty\}.\]
It is easy to check that $\mathcal{E}_1$ is a self-adjoint subspace of $\oplus_{i\in\mathbb{N}}\mathcal{V}_i$ containing $\{1\}$, and so 
is a matix order unit space, and that $\mathcal{Q}$ is a matrix Lipschitz seminorm on $\mathcal{E}_1$. 

For the evident identifications, we have
\[\mathcal{CS}(\mathcal{V}_i)\subseteq\mathcal{CS}\left(\oplus_{j=1}^i\mathcal{V}_j\right)\subseteq
\mathcal{CS}\left(\oplus_{j=1}^n\mathcal{V}_j\right)\subseteq\mathcal{CS}(\mathcal{E}_1),\ \ \ \ 0\le i\le n.\]
Given a family of graded sets $\mathbf{S}_i=(S_{i,n}), i\in I$. We denote by $\cup_{i\in I}\mathbf{S}_i$ the graded set 
$(\cup_{i\in I}S_{i,n})$. If $\mathbf{S}=(S_n)$ is a graded set in a vector space, we denote by 
$\mathrm{mco}(\mathbf{S})$ the matrix convex hull of $\mathbf{S}$. 
Let $\mathcal{Z}=(Z_n)=\mathrm{mco}\left(\cup_{i\in\mathbb{N}}\mathcal{CS}(\mathcal{V}_i\right))$, 
$\mathcal{U}=\cup_{n\in\mathbb{N}}\mathcal{CS}(\oplus_{j=1}^n\mathcal{V}_j)$ and $\mathcal{W}_n=\mathrm{mco}(\cup_{i=1}^n\mathcal{CS}(\mathcal{V}_i))$. 

\begin{proposition}\label{pro61}
$\mathcal{Z}$ and $\mathcal{U}$ are BW-dense in $\mathcal{CS}(\mathcal{E}_1)$. In particular, $\mathcal{W}_n$ is BW-dense in 
$\mathcal{CS}(\oplus_{j=1}^n\mathcal{V}_j)$ for $n\in\mathbb{N}$.
\end{proposition}

\begin{proof}
Since the matrix polar
\[\begin{array}{rcl}
Z_n^{\pi}&=&\{a\in M_n(\mathcal{E}_1): \mathrm{Re}\ll a,,
\varphi\gg\le 1_{r\times n}\mbox{  for  all  }\varphi\in Z_r, r\in\mathbb{N}\}\\
&=&\{a\in M_n(\mathcal{E}_1): \varphi_n(\mathrm{Re}(a)
\le 1_{r\times n}\mbox{  for  all  }\varphi\in Z_r, r\in\mathbb{N}\}\\
&=&\{a\in M_n(\mathcal{E}_1): \varphi_n(1_n-
\mathrm{Re}(a))\ge 0\mbox{  for  all  }\varphi\in Z_r, r\in\mathbb{N}\}\\
&=&\{a\in M_n(\mathcal{E}_1): 1_n-\mathrm{Re}(a)\ge 0\}\\
&=&\{a\in M_n(\mathcal{E}_1): \|a_{+}\|\le 1\},   
\end{array}\]
we have
\[\begin{array}{rcl}
Z_n^{\pi\pi}&=&\{f\in M_n(\mathcal{E}_1^\ast): \mathrm{Re}\ll f,a\gg\le 1_{r\times n}\\
&&\mbox{  when  }1_r-\mathrm{Re}(a)\ge 0, a\in M_r(\oplus_{i\in\mathbb{N}}\mathcal{V}_i), r\in\mathbb{N}\}.
\end{array}\]
For $f\in Z_n^{\pi\pi}, \lambda\in\mathbb{R}$ and $a=a^\ast\in M_r(\mathcal{E}_1)$, we have that 
$1_r-\mathrm{Re}(i\lambda a)\ge 0$, and so $\mathrm{Re}\ll i\lambda a, f\gg\le 1_{r\times n}$. Thus 
$$0=\mathrm{Re}\ll i a, f\gg=\mathrm{Re}(i\ll a, f\gg)=-\mathrm{Im}\ll a, f\gg,$$
that is, $\mathrm{Im}\ll a, f\gg=0$. If $a\in M_r(\mathcal{E}_1), a\ge 0$ and $\lambda\le 0$, then 
$1_r-\mathrm{Re}(\lambda a)\ge 0$ and so $\ll\lambda a, f\gg=\mathrm{Re}\ll\lambda a,f\gg\le 1_{r\times n}$. 
Thus $\ll a, f\gg\ge 0$. Clearly, $\ll 1_r,f\gg\le 1_{r\times n}$. Therefore,
$$Z_n^{\pi\pi}=\{f\in M_n(\mathcal{E}_1^\ast): f \mbox{ is completely positive and } f(1)\le 1_n\}.$$ 
By the bipolar theorem in matrix convexity (see Corollary 5.5 in \cite{efwi}), $\mathrm{mco}(\mathcal{Z}\cup\{0\})$ is BW-dense in $\mathcal{Z}^{\pi\pi}=
(Z^{\pi\pi}_n)$. Evaluting the completely positive mappings at $1$, we see that $\mathcal{Z}$ is BW-dense in $\mathcal{CS}(\mathcal{E}_1)$. 
Because $\mathcal{Z}\subseteq\mathcal{U}\subseteq\mathcal{CS}(\mathcal{E}_1)$, $\mathcal{U}$ is BW-dense in $\mathcal{CS}(\mathcal{E}_1)$.
\end{proof}

Define $\mathcal{P}_n=(P_{n,k})$ on $\oplus_{i=1}^n\mathcal{V}_i$ by
\[P_{n,k}(a_1,\cdots,a_n)=\max\{M_{i,k}(a_i,a_{i+1}): 1\le i\le n-1\},\]
for $(a_1,\cdots,a_n)\in M_k(\oplus_{i=1}^n\mathcal{V}_i)$. Similar to the proof of Proposition \ref{pro52}, we have

\begin{proposition}\label{pro62}
$\mathcal{P}_n$ is a matrix Lip-norm on $\oplus_{i=1}^n\mathcal{V}_i$, and 
induces $\mathcal{L}_j, 1\le j\le n$ and $\mathcal{M}_i$ and $\mathcal{P}_i$, $1\le i\le n-1$, via the evident projections.
\end{proposition}

For $b\in M_i(\oplus^n_{j=1}\mathcal{V}_j)$ and $\epsilon>0$, set $b_n=b$. Since $\mathcal{P}_{n+1}$ induces $\mathcal{P}_n$, we can find 
$b_{n+1}\in M_i(\oplus^{n+1}_{j=1}\mathcal{V}_j)$ such that $(\pi_n)_i(b_{n+1})=b_n$ and $P_{n+1,i}(b_{n+1})< P_{n,i}(b_n)+\frac{\epsilon}{2^n}$, where 
$\pi_n$ is the evident projection from $\oplus^{n+1}_{j=1}\mathcal{V}_j$ onto $\oplus^{n}_{j=1}\mathcal{V}_j$. 
Similarly, we can find $b_{n+2}\in M_i(\oplus^{n+2}_{j=1}\mathcal{V}_j)$ such that $(\pi_{n+1})_i(b_{n+2})=b_{n+1}$ and $P_{n+2,i}(b_{n+2})
< P_{n+1,i}(b_{n+1})+\frac{\epsilon}{2^{n+1}}$. Continuing in this way, for $t\ge n$ we get 
$b_{t+1}\in M_i(\oplus^{t+1}_{j=1}\mathcal{V}_j)$ such that $(\pi_t)_i(b_{t+1})=b_t$ and $P_{t+1,i}(b_{t+1})< P_{t,i}(b_t)+\frac{\epsilon}{2^t}$. 
We let $c=\{c_j\}$ be the unique element of $M_i(\prod_{j\in\mathbb{N}}\mathcal{V}_j)$ such that $(\pi_t)_i(c)=b_t$ for $t\ge n$. Then $Q_i(c)\le 
P_{n,i}(b)+\epsilon$. So, $P_{n,i}(b)=Q_{\oplus_{j=1}^n\mathcal{V}_j,i}(b)$. Set $d_k=(c_1,\cdots,c_k), k\in\mathbb{N}$.

Since $\mathcal{P}_n$ induces $\mathcal{M}_i(1\le i\le n-1)$, via the evidence projections, 
$$\mathrm{dist}_H^{D_{P_{n,i}}}(CS_i(\mathcal{V}_k), CS_i(\mathcal{V}_{k+1}))=\mathrm{dist}_H^{D_{M_{k,i}}}(CS_i(\mathcal{V}_k), CS_i(\mathcal{V}_{k+1})).$$ 
For $m,n\in\mathbb{N}$ with $m<n$ and 
$\varphi_n\in CS_i(\mathcal{V}_n)$, we can find $\varphi_{n-1}\in CS_i(\mathcal{V}_{n-1})$ with
\[D_{P_{n,i}}(\varphi_{n-1},\varphi_{n})\le\mathrm{dist}_H^{D_{M_{n-1,i}}}(CS_i(\mathcal{V}_{n-1}),CS_i(\mathcal{V}_n)).\]
Similarly, we can find $\varphi_{n-2}\in CS_i(\mathcal{V}_{n-2})$ with
\[D_{P_{n,i}}(\varphi_{n-2},\varphi_{n-1})\le\mathrm{dist}_H^{D_{M_{n-2,i}}}(CS_i(\mathcal{V}_{n-2}),CS_i(\mathcal{V}_{n-1})).\]
Inductively, we can find $\varphi_m,\cdots,\varphi_{n-1}$ with $\varphi_k\in CS_i(\mathcal{V}_k)$ and 
\[D_{P_{n,i}}(\varphi_{k},\varphi_{k+1})\le\mathrm{dist}_H^{D_{M_{k,i}}}(CS_i(\mathcal{V}_{k}),CS_i(\mathcal{V}_{k+1})),\]
for $m\le k\le n-1$. Consequently,
\[D_{P_{n,i}}(\varphi_{m},\varphi_{n})\le\sum_{j=m}^{n-1}\mathrm{dist}_H^{D_{M_{j,i}}}(CS_i(\mathcal{V}_{j}),CS_i(\mathcal{V}_{j+1})),\ \ \ 2\le k\le n-1.\]
Similarly, for $\varphi_m\in CS_i(\mathcal{V}_m)$ we can find a $\varphi_n\in CS_i(\mathcal{V}_n)$ such that the inequality above 
holds. Thus by Proposition 6.1, we have
\[\mathrm{dist}_H^{D_{P_{n,i}}}\left(CS_i(\oplus_{j=1}^m\mathcal{V}_j), CS_i(\mathcal{V}_n)\right)
\le\sum_{j=m}^{n-1}\mathrm{dist}_H^{D_{M_{j,i}}}(CS_i(\mathcal{V}_{j}),CS_i(\mathcal{V}_{j+1})).\]

\begin{proposition}\label{pro63}
For $m<n$, we have
\[\mathrm{dist}_H^{D_{P_{n,i}}}\left(CS_i(\oplus_{j=1}^m\mathcal{V}_j), CS_i(\oplus_{j=1}^n\mathcal{V}_j)\right)
\le\sum_{k=1}^i\sum_{j=m}^{n-1}\mathrm{dist}_H^{D_{M_{j,k}}}(CS_k(\mathcal{V}_{j}),CS_k(\mathcal{V}_{j+1})).\]
\end{proposition}

\begin{proof}
For $\varphi\in CS_i(\oplus_{j=1}^m\mathcal{V}_j)$, we can find a $\psi\in CS_i(\mathcal{V}_n)\subseteq CS_i(\oplus_{j=1}^n\mathcal{V}_j)$ such that 
$D_{P_{n,i}}(\varphi,\psi)\le\sum_{k=1}^i\sum_{j=m}^{n-1}\mathrm{dist}_H^{D_{M_{j,k}}}(CS_k(\mathcal{V}_{j}),CS_k(\mathcal{V}_{j+1}))$ from 
the discussion before the proposition. 
Suppose $\varphi\in CS_i(\oplus_{j=1}^n\mathcal{V}_j)$ and $\epsilon>0$. For each $i\in\mathbb{N}$, $\cup_{j=1}^n CS_i(\mathcal{V}_j)$ is a BW-closed 
subset of $CS_i(\oplus_{j=1}^n\mathcal{V}_j)$, and $\gamma^\ast(\cup_{j=1}^n CS_i(\mathcal{V}_j))\gamma\subseteq\cup_{j=1}^n CS_k(\mathcal{V}_j)$ for all 
isometries $\gamma\in M_{i,k}$. From Proposition \ref{pro61}, the BW-closure $\overline{\mathrm{mco}}(\cup_{j=1}^n CS_i(\mathcal{V}_j))$ of 
$\mathrm{mco}(\cup_{j=1}^n CS_i(\mathcal{V}_j))$ is $\mathcal{CS}(\oplus_{j=1}^n\mathcal{V}_j)$, and so by Theorem 4.6 and Theorem 4.3 in \cite{wewi}, there exist 
$\varphi_k\in CS_{l_k}(\mathcal{V}_{j_k})$ and $\gamma_k\in M_{l_k,i}$ for $k=1,2,\cdots,s, 1\le j_k\le n, 1\le l_k\le i$ 
satisfying $\sum_{k=1}^s\gamma^\ast_k\gamma_k=1_i$ such that
\[D_{P_{n,i}}\left(\varphi,\sum_{k=1}^s\gamma_k^\ast\varphi_k\gamma_k\right)<\epsilon.\]
For each $\varphi_k$, we can find $\psi_k\in CS_{l_k}(\oplus_{j=1}^m\mathcal{V}_j)$ so that 
$$D_{P_{n,l_k}}(\varphi_k, \psi_k)\le\sum_{j=m}^{n-1}\mathrm{dist}_H^{D_{M_{j,l_k}}}
\left(CS_{l_k}(\mathcal{V}_j),CS_{l_k}(\mathcal{V}_{j+1})\right).$$
Thus
\[\begin{array}{rcl}
&&D_{P_{n,i}}\left(\varphi, \sum_{k=1}^s\gamma^\ast_k\psi_k\gamma\right)\\
&\le&D_{P_{n,i}}\left(\varphi, \sum_{k=1}^s\gamma^\ast_k\varphi_k\gamma\right)+
D_{P_{n,i}}\left(\sum_{k=1}^s\gamma^\ast_k\varphi_k\gamma, \sum_{k=1}^s\gamma^\ast_k\psi_k\gamma\right)\\
&<&\epsilon+D_{P_{n,i}}\left(\sum_{k=1}^s\gamma^\ast_k\varphi_k\gamma, \sum_{k=1}^s\gamma^\ast_k\psi_k\gamma\right)\\
&\le&\epsilon+D_{P_{n,\sum_{k=1}^sl_k}}(\varphi_1\oplus\cdots\oplus\varphi_s, \psi_1\oplus\cdots\oplus\psi_s)\\
&=&\epsilon+\max\left\{D_{P_{n,l_1}}(\varphi_1,\psi_1),\cdots, D_{P_{n,l_s}}(\varphi_s,\psi_s)\right\}\\
&\le&\epsilon+\sum_{k=1}^i\sum_{j=m}^{n-1}\mathrm{dist}_H^{D_{M_{j,k}}}(CS_k(\mathcal{V}_{j}),CS_k(\mathcal{V}_{j+1})),
\end{array}\] 
because $\mathcal{D}_{\mathcal{P}_n}=(D_{P_{n,k}})$ is a matrix metric (see Example 5.2 in \cite{wu2}). Since $\epsilon$ is arbitrary, we obtain the desired inequality.
\end{proof}

Now for $\varphi,\psi\in CS_i(\oplus_{j=1}^n\mathcal{V}_j)$, there are $\varphi_1,\psi_1\in CS_i(\mathcal{V}_1)$ with
\[D_{P_{n,i}}(\varphi_{1},\varphi)\le\sum_{j=1}^{n-1}\mathrm{dist}_H^{D_{M_{j,i}}}(CS_i(\mathcal{V}_{j}),CS_i(\mathcal{V}_{j+1})),\]
and 
\[D_{P_{n,i}}(\psi_{1},\psi)\le\sum_{j=1}^{n-1}\mathrm{dist}_H^{D_{M_{j,i}}}(CS_i(\mathcal{V}_{j}),CS_i(\mathcal{V}_{j+1})).\]
So
\[\begin{array}{rcl}
D_{P_{n,i}}(\varphi,\psi)&\le&D_{P_{n,i}}(\varphi,\varphi_1)+D_{P_{n,i}}(\varphi_1,\psi_1)+D_{P_{n,i}}(\psi_1,\psi)\\
&\le&\mathrm{diam}(CS_i(\mathcal{V}_1),D_{L_{1,i}})+2\sum_{j=1}^{n-1}\mathrm{dist}_H^{D_{M_{j,i}}}(CS_i(\mathcal{V}_{j}),CS_i(\mathcal{V}_{j+1}))\\
&\buildrel \triangle\over =&h_i,
\end{array}\]
where $\mathrm{diam}(CS_i(\mathcal{V}_1),D_{L_{1,i}})$ is the diameter of $CS_i(\mathcal{V}_1)$ with respect to $D_{L_{1,i}}$.

By Proposition 5.2 and Proposition 3.8 in \cite{wu1}, we have 
\[\|\tilde{d}_n\|^\sim_i\le h_iP_{n,i}(d_n)\le h_iQ_i(c)< h_iQ_i(c)+\epsilon,\]
where $\epsilon>0$. So there is an $\alpha_{n,i}=\left[\alpha_{st}^{(n,i)}\right]\in M_i$ such that
\[\left\|d_n-\left[\alpha_{st}^{(n,i)}(1,\cdots,1)\right]\right\|_i\le h_iQ_i(c)+\epsilon,\ \ \ \ n\in\mathbb{N}.\]
Set
\[G_{n,i}=\left\{\beta_{n,i}=\left[\beta_{st}^{(n,i)}\right]\in M_i: \left\|d_n-\left[\beta_{st}^{(n,i)}(1,\cdots,1)\right]\right\|_i\le h_iQ_i(c)+\epsilon\right\}.\]
Then $G_{n,i}$ is a non-empty closed bounded subset of $M_i$. Clearly, $G_{n+1,i}\subseteq G_{n,i}$. So there exists a $\beta_0\in
\cap_{n=1}^\infty G_{n,i}$. We have
\[\|d_n\|_i\le\|\beta_0\|+h_iQ_i(c)+\epsilon,\ \ \ \ n\in\mathbb{N}.\]
Thus $c\in M_i(\oplus_{j\in\mathbb{N}}\mathcal{V}_j)$, and we obtain

\begin{proposition}\label{pro64}
For $n\in\mathbb{N}$, $\mathcal{Q}$ induces $\mathcal{P}_n$ via the evident projection.
\end{proposition}

\begin{theorem}\label{th65}
The metric space $\mathfrak{R}$ of complete isometry classes of quantized metric spaces, with the metric $\mathrm{dist}_{NC}$, is 
complete.
\end{theorem}

\begin{proof}
Let $\{(\mathcal{V}_n, \mathcal{L}_n)\}$ be a sequence in $\mathfrak{R}$ which is Cauchy with respect to the quantized 
Gromov-Hausdorff distance $\mathrm{dist}_{NC}$. To show that $\{(\mathcal{V}_n, \mathcal{L}_n)\}$ converges it suffices to show 
that a subsequence converges. Since $\{(\mathcal{V}_n, \mathcal{L}_n)\}$ is Cauchy, we can choose a subsequence, still denoted by 
$\{(\mathcal{V}_n, \mathcal{L}_n)\}$, such that 
$$\mathrm{dist}_{NC}(\mathcal{V}_n,\mathcal{V}_{n+1})<\frac1{2^n},$$
for all $n\in\mathbb{N}$. By definition, there exist $\mathcal{M}_n=(M_{n,k})\in\mathcal{M}(\mathcal{L}_{n}, \mathcal{L}_{n+1})$ with
\[\sup_{k\in\mathbb{N}}\left\{k^{-2}\mathrm{dist}_H^{D_{M_{n,k}}}(CS_k(\mathcal{V}_{n}), CS_k(\mathcal{V}_{n+1}))\right\}<\frac1{2^n},\]
for all $n\in\mathbb{N}$. It follows that
\[\sum_{n=1}^\infty \sup_{k\in\mathbb{N}}\left\{k^{-2}\mathrm{dist}_H^{D_{M_{n,k}}}(CS_k(\mathcal{V}_{n}), CS_k(\mathcal{V}_{n+1}))\right\}<+\infty.\]

Let $\epsilon>0$ be given. Then there is an $m\in\mathbb{N}$ such that
\[\sum_{n=m}^\infty\sup_{k\in\mathbb{N}}\left\{k^{-2}\mathrm{dist}_H^{D_{M_{n,k}}}(CS_k(\mathcal{V}_{n}), CS_k(\mathcal{V}_{n+1}))\right\}<\epsilon.\]
By Proposition \ref{pro62}, Proposition \ref{pro63} and Proposition \ref{pro64}, we have
\[\begin{array}{rcl}
&&\mathrm{dist}_H^{D_{Q_i}}\left(CS_i(\oplus_{j=1}^m\mathcal{V}_{j}), CS_i(\oplus_{j=1}^n\mathcal{V}_{j})\right)\\
&\le&\sum_{k=1}^i\sum_{j=m}^{n-1}\mathrm{dist}_H^{D_{M_{j,k}}}(CS_k(\mathcal{V}_{j}),CS_k(\mathcal{V}_{j+1}))\\
&\le&\sum_{k=1}^ik^2\sum_{j=m}^{n-1}\sup_{k\in\mathbb{N}}\left\{k^{-2}\mathrm{dist}_H^{D_{M_{j,k}}}(CS_k(\mathcal{V}_{j}),CS_k(\mathcal{V}_{j+1}))\right\}\\
&<&(\sum^i_{k=1}k^2)\epsilon,\end{array}\]
for $n>m$. This says that $CS_i(\oplus_{j=1}^m\mathcal{V}_{j})$ is $(\sum^i_{k=1}k^2)\epsilon$-dense for $D_{Q_i}$ in $Z_i$. But 
$CS_i(\oplus_{j=1}^m\mathcal{V}_{j})$ is BW-compact for the topology from $D_{Q_i}=D_{P_{m,i}}$ by Proposition \ref{pro62}. Thus 
$CS_i(\oplus_{j=1}^m\mathcal{V}_{j})$ is totally bounded for $D_{Q_i}$, and so $Z_i$ is totally bounded for $D_{Q_i}$. Let $\hat{\mathcal{Z}}=(\hat{Z}_n)$ 
be the completion of $\mathcal{Z}$ for $\mathcal{D}_{\mathcal{Q}}$. We let $\mathcal{D}_{\mathcal{Q}}$ denote also the extension of 
$\mathcal{D}_{\mathcal{Q}}$ to $\hat{\mathcal{Z}}$. Then $\hat{\mathcal{Z}}$ is a compact matrix convex set.

For $\{a_i\}\in M_n(\mathcal{E}_1)$, we have
\[\ll\sum_{j=1}^m\gamma^\ast_j\varphi_j\gamma_j,\{a_i\}\gg=\sum_{j=1}^m(\gamma_j\otimes 1_n)^\ast\ll\varphi_j,\{a_i\}\gg(\gamma_j\otimes 1_n),\]
and 
\[\begin{array}{rcl}
&&\left\|\ll\sum_{j=1}^m\gamma^\ast_j\varphi_j\gamma_j,\{a_i\}\gg-\ll\sum_{k=1}^p\lambda^\ast_k\psi_k\lambda_k,\{a_i\}\gg\right\|\\
&\le&L_{D_{Q_n}}(\{a_i\})D_{Q_r}\left(\sum_{j=1}^m\gamma^\ast_j\varphi_j\gamma_j,\sum_{k=1}^p\lambda^\ast_k\psi_k\lambda_k\right)\\
&\le&Q_n(\{a_i\})D_{Q_r}\left(\sum_{j=1}^m\gamma^\ast_j\varphi_j\gamma_j,\sum_{k=1}^p\lambda^\ast_k\psi_k\lambda_k\right),
\end{array}\]
where $\varphi_j\in CS_{n_j}(\mathcal{V}_{q_j})$, $\psi_k\in CS_{m_k}(\mathcal{V}_{l_k})$, and $\gamma_j\in M_{n_j,r}$, and $\lambda_k\in M_{m_k,r}$
satisfying $\sum_{j=1}^m\gamma^\ast_j\gamma_j=1_r$ and $\sum_{k=1}^p\lambda^\ast_k\lambda_k=1_r$. So the map $\Phi: \mathcal{E}_1\mapsto A(\mathcal{Z}
)$, given by
\[(\Phi(\{a_i\}))\left(\sum_{j=1}^m\gamma^\ast_j\varphi_j\gamma_j\right)=\sum_{j=1}^m\gamma^\ast_j\varphi_j(a_{q_j})\gamma_j,\]
for $\{a_i\}\in\mathcal{E}_1$, $\varphi_j\in CS_{n_j}(\mathcal{V}_{q_j})$ and $\gamma_j\in M_{n_j,r}$ satisfying 
$\sum_{j=1}^m\gamma^\ast_j\gamma_j=1_r$, is well-defined and $\Phi(\{a_i\})$ can be extended to an element $\widehat{\Phi(\{a_i\})}\in A(\hat{
\mathcal{Z}})$. Moreover if $\{a_i\}\ge 0$ in $\mathcal{E}_1$ then $\widehat{\Phi(\{a_i\})}\ge 0$ in $A(\hat{\mathcal{Z}})$ and $\widehat{\Phi(\{1\})}
=\mathbf{I}$. Thus $\mathcal{E}_1$ can be regarded as a matrix order unit subspace of $A(\hat{\mathcal{Z}})$.

Define the map $\Psi_r: \hat{Z}_r\mapsto CS_r(\mathcal{E}_1), r\in\mathbb{N}$, by
\[\Psi_r(z)(\{a_i\})=\widehat{\Psi(\{a_i\})}(z),\]
for $z\in\hat{Z}_r$ and $\{a_i\}\in\mathcal{E}_1$. Clearly, $\Psi$ is continuous. For $z=\sum_{j=1}^m\gamma^\ast_j\varphi_j\gamma_j\in Z_r$ with 
$\varphi_j\in CS_{n_j}(\mathcal{V}_{q_j})$, $\gamma_j\in M_{n_j,r}$ satisfying $\sum_{j=1}^m\gamma^\ast_j\gamma_j=1_r$, we have
\[\Psi_r(z)(\{a_i\})=\widehat{\Psi(\{a_i\})}(z)=\Psi(\{a_i\})(z)=z(\{a_i\}),\]
that is, $\Psi_r(z)=z$. Since $Z_r$ is dense in $CS_r(\mathcal{E}_1)$ and $\hat{Z}_r$ is compact, we obtain that $\Psi_r(\hat{Z}_r)=CS_r(\mathcal{E}_1)$.

If $z_1, z_2\in\hat{Z}_r$ with $z_1\neq z_2$ and $k=D_{Q_r}(z_1,z_2)$, we can find $y_1, y_2\in Z_r$ such that $D_{Q_r}(z_i,y_i)<\frac k4, i=1,2$. 
Thus $D_{Q_r}(y_1,y_2)>\frac k2$. So we can find $\{w_i\}\in M_r(\mathcal{E}_1)$ with $Q_r(\{w_i\})\le 1$ and $\|\ll\{w_i\},y_1\gg-\ll\{w_i\},y_2\gg\|
>\frac k2$. But $L_{D_{Q_r}}(\widehat{\Phi_r(\{w_i\})}\le 1$ so that $\|\ll\widehat{\Phi_r(\{w_i\})},z_i\gg-\ll\widehat{\Phi_r(\{w_i\})},y_i\gg\|<
\frac k4, i=1,2$. Thus $\|\ll\widehat{\Phi_r(\{w_i\})},z_1\gg-\ll\widehat{\Phi_r(\{w_i\})},z_2\gg\|>0$. Therefore, $\Psi_r$ is injective. So $\Psi_r$ is a 
homeomorphism of $\hat{Z}_r$ onto $CS_r(\mathcal{E}_1)$ for $r\in\mathbb{N}$. From this we see that the $\mathcal{D}_{\mathcal{Q}}$-topology on 
$\mathcal{CS}(\mathcal{E}_1)$ agrees with the BW-topology. Hence $\mathcal{Q}$ is a matrix Lip-norm on $\mathcal{E}_1$.

By Proposition \ref{pro62} and Proposition \ref{pro64}, we obtain
\[\sum_{n=1}^\infty\sup_{k\in\mathbb{N}}\left\{k^{-2}\mathrm{dist}_H^{D_{Q_k}}(CS_k(\mathcal{V}_n), CS_k(\mathcal{V}_{n+1}))\right\}<+\infty,\]
which indicate that, for $k\in\mathbb{N}$, $\{CS_k(\mathcal{V}_n)\}$ is a Cauchy sequence for $\mathrm{dist}^{D_{Q_k}}_H$, and has a limit $K_k\subseteq CS_k(\mathcal{E}_1)$. Clearly 
$\mathcal{K}=(K_k)$ is a compact matrix convex set.

Because $\mathcal{E}_1$ is completely order isomorphic to a dense subspce of $A(\mathcal{CS}(\mathcal{E}_1))$ (Proposition 6.1(1) in \cite{wu2}), 
we can view $\mathcal{E}_1$ as a dence subspace of $A(\mathcal{CS}(\mathcal{E}_1))$. Let $\phi$ be the map which restricts the elements of $A(\mathcal{CS}(\mathcal{E}_1))$ 
to $\mathcal{K}$ and $\mathcal{V}=\phi(\mathcal{E}_1)$. Then $(\mathcal{V},\mathcal{Q}_{\mathcal{V}})$ is a quantized metric space. 

Given $\epsilon>0$. Then there is an $N\in\mathbb{N}$ such that
\[\sum_{j=n}^\infty\sup_{k\in\mathbb{N}}\left\{k^{-2}\mathrm{dist}_H^{D_{Q_k}}(CS_k(\mathcal{V}_j), CS_k(\mathcal{V}_{j+1}))\right\}<\epsilon, \ \ \ n\ge N.\]
For $k,p\in\mathbb{N}$, we have
\[\begin{array}{rcl}
&&k^{-2}\mathrm{dist}_H^{D_{Q_k}}(CS_k(\mathcal{V}_n), CS_k(\mathcal{V}_{n+p}))\\
&\le&\sum_{j=n}^{n+p-1}k^{-2}\mathrm{dist}_H^{D_{Q_k}}(CS_k(\mathcal{V}_j), CS_k(\mathcal{V}_{j+1}))\\
&\le&\sum_{j=n}^\infty\sup_{k\in\mathbb{N}}\left\{k^{-2}\mathrm{dist}_H^{D_{Q_k}}(CS_k(\mathcal{V}_j), CS_k(\mathcal{V}_{j+1}))\right\}\\
&<&\epsilon,\end{array}\]
for $n\ge\mathbb{N}$. Letting $p\to +\infty$, we obtain
\[k^{-2}\mathrm{dist}_H^{D_{Q_k}}(CS_k(\mathcal{V}_n), K_k)\le\epsilon,\]
for $k\in\mathbb{N}$, and so $\sup_{k\in\mathbb{N}}\left\{k^{-2}\mathrm{dist}_H^{D_{Q_k}}(CS_k(\mathcal{V}_n), K_k)\right\}\le\epsilon$. 
By Proposition \ref{pro475}, for $n\ge N$ we have 
\[\mathrm{dist}_{NC}(\mathcal{V}_n,\mathcal{V})\le\sup_{k\in\mathbb{N}}\left\{k^{-2}\mathrm{dist}_H^{D_{Q_k}}(CS_k(\mathcal{V}_n), 
K_k)\right\}\le\epsilon.\]
Therefore, $\lim_{n\to\infty}\mathrm{dist}_{NC}(\mathcal{V}_n,\mathcal{V})=0$.
\end{proof}

\section{Matrix approximability}\label{s7}

In this section, we establish a matrix approximability theorem for $1$-exact matrix order unit spaces.

\begin{lemma}\label{le701}
Let $(\mathcal{V}, \mathcal{L})$ be a quantized metric space and let $x=[x_{st}]\in M_k(\mathcal{V})$, $x_{st}=x_{st}^{(1)}+i
x_{st}^{(2)}$ with $(x_{st}^{(p)})^\ast=x_{st}^{(p)}$ for $p=1,2, s,t=1,2,\cdots,k$. Suppose $\lambda_{st}^{(p)}\in\sigma(x_{st}^{(p)})$ 
for $p=1,2, s,t=1,2,\cdots,k$. Then $\|x-[(\lambda_{st}^{(1)}+i\lambda_{st}^{(2)})1]\|_k\le 2k^2L_k(x)\mathrm{diam}(\mathcal{V},\mathcal{L})$.
\end{lemma}

\begin{proof}
By Proposition 2.11 in \cite{ker}, we have
\[\begin{array}{rcl}
\left\|x-\left[\left(\lambda_{st}^{(1)}+i\lambda_{st}^{(2)}\right)1\right]\right\|_k
&\le&\left\|\left[x_{st}^{(1)}-\lambda_{st}^{(1)}1\right]\right\|_k+\left\|\left[x_{st}^{(2)}-
\lambda_{st}^{(2)}1\right]\right\|_k\\
&\le&\sum_{p=1}^2\sum_{s,t=1}^k\left\|x_{st}^{(p)}-\lambda_{st}^{(p)}1\right\|_1\\
&\le&\sum_{p=1}^2\sum_{s,t=1}^k L_1(x_{st}^{(p)})\mathrm{diam}(\mathcal{V}, \mathcal{L})\\
&\le&\sum_{s,t=1}^k2L_k(x)\mathrm{diam}(\mathcal{V}, \mathcal{L})\\
&=&2k^2L_k(x)\mathrm{diam}(\mathcal{V},\mathcal{L}).
\end{array}\]
\end{proof}

An operator space $\mathcal{X}$ is said to be $1$-{\it exact} if for every finite-dimensional subspace $\mathcal{E}\subseteq\mathcal{X}$ and 
$\lambda>1$ there is an isomorphism $\alpha$ from $\mathcal{E}$ onto a subspace of a matrix algebra such that 
$\|\alpha\|_{cb}\|\alpha^{-1}\|_{cb}\le\lambda$. A matrix order unit space $(\mathcal{V},1)$ is said to be $1$-{\it exact} if it is 
$1$-exact as an operator space.

\begin{theorem}\label{pro73}
Let $(\mathcal{V}, \mathcal{L})$ be a quantized metric space. If $\mathcal{V}$ is $1$-{\it exact}, then for every $\epsilon>0$, there 
is a quantized metric space $(M_{n_{\lambda_\epsilon}}, \mathcal{N})$ such that 
$$\mathrm{dist}_{NC}(\mathcal{V}, M_{n_{\lambda_\epsilon}})<\epsilon.$$
\end{theorem}

\begin{proof}
Since $\mathcal{V}$ is $1$-exact, by Lemma 5.1 in \cite{keli} there is a unital complete order embedding $\iota: \mathcal{V}\mapsto
\mathscr{B}(\mathscr{H})$ and a net
\[\mathcal{V}\xrightarrow{\varphi_\lambda}M_{n_{\lambda}}\xrightarrow{\psi_\lambda}\mathscr{B}(\mathscr{H})\]
of unital completely positive mappings through matrix algebras such that $\psi_\lambda\circ\varphi_\lambda$ converges pointwise to $\iota$. 
Given $\epsilon>0$. By Lemma 7.2, we have
\[L_1^1=B_1^{2\mathrm{diam}(\mathcal{V},\mathcal{L})}+\mathbb{C}1,\]
where $B_1^{2\mathrm{diam}(\mathcal{V},\mathcal{L})}=\{a\in\mathcal{V}: L_1(a)\le 1, \|a\|_1\le 2\mathrm{diam}(\mathcal{V}, \mathcal{L})\}$. 
From Proposition 7.5 in \cite{wu2}, $B_1^{2\mathrm{diam}(\mathcal{V},\mathcal{L})}$ is totally bounded for $\|\cdot\|_1$. So there is a 
$\lambda_{\epsilon}$ such that
\[\|(\psi_{\lambda_{\epsilon}}\circ\varphi_{\lambda_{\epsilon}})(x)-x\|<\frac{\epsilon}5,\ \ \ \ x\in L_1^1.\]
Denote $\mathcal{W}=\varphi_{\lambda_{\epsilon}}(\mathcal{V})$ and $Q_k(y)=\inf\{L_k(x): (\varphi_{\lambda_{\epsilon}})_k(x)=y\}$ 
for $y\in M_k(\mathcal{W})$ and $k\in\mathbb{N}$.

We define
\[N_k(x,y)=\frac5{\epsilon}\|(\varphi_{\lambda_\epsilon})_k(x)-y\|_k,\ \ \ (x,y)\in M_k(\mathcal{V}\oplus\mathcal{W}), k\in\mathbb{N}.\]
It is clear that $\mathcal{N}=(N_k)$ is a matrix seminorm on $\mathcal{V}\oplus\mathcal{W}$ and satisfies the conditions (1), (2) and (3) of 
Definition \ref{de43}. For $x\in M_k(\mathcal{V})$ and $\delta>0$, we can choose $y=(\varphi_{\lambda_\epsilon})_k(x)$. Then
\[\max\{Q_k(y), N_k(x,y)\}=Q_k(y)\le L_k(x)\le L_k(x)+\delta.\]
For $y\in M_k(\mathcal{W})$ and $\delta>0$, we can take $x\in M_k(\mathcal{V})$ such that $y=(\varphi_{\lambda_\epsilon})_k(x)$ and $L_k(x)\le Q_k(y)+\delta$. 
Then
\[\max\{L_k(x), N_k(x,y)\}=L_k(x)\le Q_k(y)+\delta.\]
So $\mathcal{N}$ is a matrix bridge between  $(\mathcal{V}, \mathcal{L})$ and $(\mathcal{W}, \mathcal{Q})$. Define 
\[P_k(x,y)=\max\{L_k(x), Q_k(y), N_k(x,y)\},\ \ \ (x,y)\in M_k(\mathcal{V}\oplus\mathcal{W}), k\in\mathbb{N}.\]
Then $\mathcal{P}=(P_k)\in \mathcal{M}(\mathcal{L},\mathcal{Q})$ by Proposition \ref{pro52}. If $f\in CS_k(\mathcal{W})$, we have 
$f\circ\varphi_{\lambda_\epsilon}\in CS_k(\mathcal{V})$ and
\[\begin{array}{rcl}
&&D_{P_k}(f,f\circ\varphi_{\lambda_\epsilon})\\
&=&\sup\{\|\ll f,y\gg-\ll f\circ\varphi_{\lambda_\epsilon},x\gg\|: P_k(x,y)\le 1, (x,y)\in M_k(\mathcal{V}\oplus\mathcal{W})\}\\
&=&\sup\{\|\ll f,y-(\varphi_{\lambda_\epsilon})_k(x)\gg\|: P_k(x,y)\le 1, (x,y)\in M_k(\mathcal{V}\oplus\mathcal{W})\}\\
&\le&\sup\{\|y-(\varphi_{\lambda_\epsilon})_k(x)\|_k: P_k(x,y)\le 1, (x,y)\in M_k(\mathcal{V}\oplus\mathcal{W})\}\\
&\le&\frac{\epsilon}5
\end{array}\]
On the other hand, if $g\in CS_k(\mathcal{V})$, $g$ can be extended to a $\overline{g}\in CS_k(\mathscr{B}(\mathscr{H}))$ by Arveson's 
extension theorem. We have $\overline{g}\circ\psi_{\lambda_\epsilon}\in CS_k(\mathcal{W})$ and 
\[\begin{array}{rcl}
&&D_{P_k}(g,\overline{g}\circ\psi_{\lambda_\epsilon})\\
&=&\sup\{\|\ll g,x\gg-\ll\overline{g}\circ\psi_{\lambda_\epsilon},y\gg\|: P_k(x,y)\le 1, (x,y)\in M_k(\mathcal{V}\oplus\mathcal{W})\}\\
&=&\sup\{\|\ll\overline{g},x-(\psi_{\lambda_\epsilon})_k(y)\gg\|: P_k(x,y)\le 1, (x,y)\in M_k(\mathcal{V}\oplus\mathcal{W})\}\\
&\le&\sup\{\|\ll\overline{g},x-(\psi_{\lambda_\epsilon}\circ\varphi_{\lambda_\epsilon})_k(x)\gg\|
+\|\ll\overline{g}\circ\psi_{\lambda_\epsilon}, (\varphi_{\lambda_\epsilon})_k(x)-y\gg\|:\\ 
&&P_k(x,y)\le 1, (x,y)\in M_k(\mathcal{V}\oplus\mathcal{W})\}\\
&\le&\sup\{\|x-(\psi_{\lambda_\epsilon}\circ\varphi_{\lambda_\epsilon})_k(x)\|_k
+\|(\varphi_{\lambda_\epsilon})_k(x)-y\|_k:\\ 
&&P_k(x,y)\le 1, (x,y)\in M_k(\mathcal{V}\oplus\mathcal{W})\}\\
&\le&\sup\big\{\sum_{i,j=1}^k\|x_{ij}-(\psi_{\lambda_\epsilon}\circ\varphi_{\lambda_\epsilon})(x_{ij})\|_1: 
P_k(x,y)\le 1,\\
&&(x,y)\in M_k(\mathcal{V}\oplus\mathcal{W})\big\}+\frac{\epsilon}5\\
&\le&\left(k^2+1\right)\frac{\epsilon}5
\end{array}\]
So we obtain that $\mathrm{dist}_{NC}(\mathcal{V},\mathcal{W})<\frac{\epsilon}2$.

Since $\mathcal{W}\subseteq M_{n_{\lambda_\epsilon}}$ is finite-dimensional, $K=Q^1_1$ is a normed-closed (and hence 
weakly closed)  absolutely convex set in $M_{n_{\lambda_\epsilon}}$, and $\mathcal{Q}^1=(Q^1_k)$ is a normed-closed (and hence 
weakly closed) absolutely matrix convex set in $M_{n_{\lambda_\epsilon}}$. Then for the corresponding matrix seminorm 
$\check{\mathcal{R}}=(\check{R}_k)$ of the maximal envelope $\check{\mathcal{K}}$ of $K$ in $M_{n_{\lambda_\epsilon}}$(see Example 
\ref{ex33}), we have 
\[\check{R}^1_1=Q^1_1,\ \ \left.\check{R}_k\right|_{M_k(\mathcal{W})}\le\left.Q_k\right|_{M_k(\mathcal{W})},\ \ k\in\mathbb{N},\]
(see page 181 in \cite{efwe}). It is clear that $\check{\mathcal{R}}$ is a matrix Lipschitz seminorm. Since the image of 
$Q^1_1=\check{R}^1_1$ in $\mathcal{W}/(\mathbb{C}1)$ is totally bounded for $\|\cdot\|^\sim$ and $\mathcal{W}\subseteq 
M_{n_{\lambda_\epsilon}}$, the image of $\check{R}^1_1$ in $M_{n_{\lambda_\epsilon}}/(\mathbb{C}1)$ is totally bounded for $\|\cdot\|^\sim$. 
By Theorem 5.3 in \cite{wu1}, $\mathcal{D}_{\mathcal{R}}$-topology on $\mathcal{CS}(M_{n_{\lambda_\epsilon}})$ agrees with the 
BW-topology. So $\check{\mathcal{R}}$ is a matrix Lip-norm on $(M_{n_{\lambda_\epsilon}},1)$. By Lemma 3.2.3 in \cite{blac}, there is a (real linear) 
projection $T$ from $(M_{n_{\lambda_\epsilon}})_{sa}$ onto $(\mathcal{W})_{sa}$ with $\|T\|\le n_{\lambda_\epsilon}$. We define 
$S: M_{n_{\lambda_\epsilon}}\mapsto\mathcal{W}$ by $S(a+ib)=T(a)+iT(b)$ for $a,b\in (M_{n_{\lambda_\epsilon}})_{sa}$. Then $S$ is a bounded 
linear mapping with $\|S\|\le 2n_{\lambda_{\epsilon}}$. Define 
\[N_k(x)=\max\left\{Q_k(S_k(x)), \check{R}_k(x), \frac4{\epsilon}\|x-S_k(x)\|_k\right\},\ \ \ x\in M_k(M_{n_{\lambda_\epsilon}}), k\in\mathbb{N}.\]
It is clear that $\mathcal{N}=(N_k)$ is a matrix Lip-norm on $M_{n_{\lambda_\epsilon}}$ since $\check{R}_k\le N_k$ for all $k\in\mathbb{N}$ 
and $\check{\mathcal{R}}$ is a matrix Lip-norm. And for $y\in M_k(\mathcal{W})$, we have
\[\begin{array}{rcl}
N_k(y)&=&\max\left\{Q_k(S_k(y)), \check{R}_k(y), \frac4{\epsilon}\|y-S_k(y)\|_k\right\}\\
&=&\max\left\{Q_k(y), \check{R}_k(y)\right\}\\
&=&Q_k(y).
\end{array}\]
Define
\[X_k(x,y)=\frac4{\epsilon}\|y-S_k(x)\|_k,\ \ \ (x,y)\in M_k(M_{n_{\lambda_\epsilon}}\oplus\mathcal{W}), k\in\mathbb{N}.\]
It is clear that $\mathcal{N}=(N_k)$ is a matrix seminorm on $M_{n_{\lambda_\epsilon}}\oplus\mathcal{W}$ and satisfies the conditions (1), 
(2) and (3) of 
Definition \ref{de43}. For $x\in M_k(M_{n_{\lambda_\epsilon}})$ and $\delta>0$, we choose $y=S_k(x)$. Then we have
\[\max\{Q_k(y), X_k(x,y)\}=Q_k(S_k(x))\le N_k(x)\le N_k(x)+\delta.\]
For $y\in M_k(\mathcal{W})$ and $\delta>0$, we choose $x=y$. Then we have
\[\max\{N_k(x), X_k(x,y)\}=N_k(y)=Q_k(y)\le Q_k(x)+\delta.\]
So $\mathcal{X}=(X_k)$ is a matrix bridge between $(M_k(M_{n_{\lambda_\epsilon}}), \mathcal{N})$ and $(\mathcal{W}, \mathcal{Q})$. 
Define 
\[Y_k(x,y)=\max\{N_k(x), Q_k(y), X_k(x,y)\},\ \ \ (x,y)\in M_k(M_{n_{\lambda_\epsilon}}\oplus\mathcal{W}), k\in\mathbb{N}.\]
Then $\mathcal{Y}=(Y_k)\in \mathcal{M}(\mathcal{N},\mathcal{Q})$ by Proposition \ref{pro52}. For 
$\varphi\in CS_k(M_{n_{\lambda_\epsilon}})$, $\psi=\varphi|_{\mathcal{W}}\in CS_k(\mathcal{W})$ and 
\[\begin{array}{rcl}
&&D_{Y_k}(\varphi,\psi)\\
&=&\sup\{\|\ll\varphi,x\gg-\ll\psi,y\gg\|: (x,y)\in M_k(M_{n_{\lambda_\epsilon}}\oplus\mathcal{W}), Y_k(x,y)\le 1\}\\
&\le&\sup\{\|\ll\varphi,x\gg-\ll\varphi,S_k(x)\gg+\ll\varphi,S_k(x)\gg-\ll\varphi,y\gg\|: \\
&&(x,y)\in M_k(M_{n_{\lambda_\epsilon}}\oplus\mathcal{W}), Y_k(x,y)\le 1\}\\
&\le&\sup\{\|x-S_k(x)\|_k+\|S_k(x)-y\|_k: (x,y)\in M_k(M_{n_{\lambda_\epsilon}}\oplus\mathcal{W}), Y_k(x,y)\le 1\}\\
&\le&\frac{\epsilon}2.
\end{array}\]
For $\psi\in CS_k(\mathcal{W})$, $\psi$ can be extended to a $\varphi\in CS_k(M_{n_{\lambda_\epsilon}})$ by Arveson's 
extension theorem. We have  
\[\begin{array}{rcl}
&&D_{Y_k}(\varphi,\psi)\\
&=&\sup\{\|\ll\varphi,x\gg-\ll\psi,y\gg\|: (x,y)\in M_k(M_{n_{\lambda_\epsilon}}\oplus\mathcal{W}), Y_k(x,y)\le 1\}\\
&\le&\sup\{\|\ll\varphi,x\gg-\ll\varphi,S_k(x)\gg+\ll\varphi,S_k(x)\gg-\ll\varphi,y\gg\|: \\
&&(x,y)\in M_k(M_{n_{\lambda_\epsilon}}\oplus\mathcal{W}), Y_k(x,y)\le 1\}\\
&\le&\sup\{\|x-S_k(x)\|+\|S_k(x)-y\|: (x,y)\in M_k(M_{n_{\lambda_\epsilon}}\oplus\mathcal{W}), Y_k(x,y)\le 1\}\\
&\le&\frac{\epsilon}2.
\end{array}\]
So $\mathrm{dist}_{NC}(\mathcal{W}, M_{n_{\lambda_\epsilon}})<\frac{\epsilon}2$. Therefore, 
$$\mathrm{dist}_{NC}(\mathcal{V}, M_{n_{\lambda_\epsilon}})\le 
\mathrm{dist}_{NC}(\mathcal{V}, \mathcal{W})
+\mathrm{dist}_{NC}(\mathcal{W}, M_{n_{\lambda_\epsilon}})<\epsilon.$$
\end{proof}

\section{Sphere as the limit of matrix algebras}\label{s8}

Let $G$ be a connected compact semisimple Lie group with a continuous length function $l$ on $G$, which satisfies the additional condition 
$l(xyx^{-1})=l(y)$ for all $x,y\in G$. Fix an irreducible unitary representation $(U,\mathcal{H})$ of $G$. 
Then $(U,\mathcal{H})$ have a highest weight vector $\xi$, of norm $1$, unique up to a scalar multiple. Let $P$ be the rank-one 
projection for $\xi$. Denote by $H$ the stability subgroup for $P$ under $\alpha$. 
For any $n\in\mathbb{N}$, we form the $n^{th}$ inner tensor power $(U^{\otimes n},\mathcal{H}^{\otimes n})$ of $(U,\mathcal{H})$. 
Let $(U^{(n)},\mathcal{H}^{(n)})$ denote the subrepresentation generated by $\xi^{(n)}=\xi^{\otimes n}$. Then 
$(U^{(n)},\mathcal{H}^{(n)})$ 
is irreducible with $\xi^{(n)}$ as highest weight vector. 
We let $\mathcal{B}^{(n)}=\mathcal{B}(\mathcal{H}^{(n)})$. The action of $G$ on $\mathcal{B}^{(n)}$ by conjugation by 
$U^{(n)}$ is denoted by $\alpha^{(n)}$.  We let $\lambda$ denote the action of $G$ on $G/H$, and so on $\mathcal{A}=C(G/H)$, by left-translation. 
We denote the corresponding Lip-norm for $\alpha^{(n)}$ and $l$ on $\mathcal{B}^{(n)}$ by $L^{(n)}$, that is,
\[L^{(n)}(T)=\sup\left\{\frac{\|\alpha^{(n)}_x(T)-T\|}{l(x)}:\ x\neq e, x\in G\right\},\ \ \ \ T\in\mathcal{B}^{(n)},\]
and we denote the Lip-norm for $\lambda$ and $l$ on $\mathcal{A}$ by $L$, that is, 
\[L(f)=\sup\left\{\frac{\|\lambda_x(f)-f\|_{\infty}}{l(x)}:\ x\neq e, x\in G\right\},\ \ \ \ f\in \mathcal{A},\]
here we view $C(G/H)$ as a subalgebra of $C(G)$. By Theorem 3.2 in \cite{ri5}, the quantum metric spaces $(\mathcal{B}^{(n)}, L^{(n)})$ 
converge to $(\mathcal{A},L)$ for quantum Gromov-Hausdorff distance as $n$ goes to $\infty$. In this section, a more general statement is established.

Let $\|\cdot\|_\infty=(\|\cdot\|_{\infty,k})$ be the matrix norm on $\mathcal{A}$. Set $\mathcal{L}^{(n)}=(L_k^{(n)})$, where
\[L_k^{(n)}(T)=\sup\left\{\frac{\|[\alpha^{(n)}_x(T_{ij})-T_{ij}]\|}{l(x)}:\ x\neq e, x\in G\right\},\ T=[T_{ij}]\in 
M_k(\mathcal{B}^{(n)}), k\in\mathbb{N},\]
and $\mathcal{L}=(L_k)$, where
\[L_k(f)=\sup\left\{\frac{\|[\lambda_x(f_{ij})-f_{ij}]\|_{\infty,k}}{l(x)}:\ x\neq e, x\in G\right\},\ f=[f_{ij}]\in 
M_k(\mathcal{A}), k\in\mathbb{N}.\]
Then $(\mathcal{B}^{(n)},\mathcal{L}^{(n)})$ and $(\mathcal{A},\mathcal{L})$ are quantized metric spaces for all $n\in\mathbb{N}$ by 
Example 6.5 in \cite{wu2}. As in \cite{ri5}, we will not restrict $\mathcal{L}$ to the Lipschitz functions. 
Let $P^{(n)}$ denote the rank-one projection for $\xi^{(n)}$. We denote the corresponding Berezin symbol mapping from $\mathcal{B}^{(n)}$ to 
$\mathcal{A}$ by $\sigma^{(n)}$. Then $\sigma^{(n)}$ is unital, positive, norm-nonincreasing and $\alpha^{(n)}$-$\lambda$-equivariant (see 
page 73 in \cite{ri5}). 
For $k\in\mathbb{N}$ and $T=[T_{ij}]\in M_k(\mathcal{B}^{(n)})$, define
\[\sigma^{(n)}_T(x)=[\sigma^{(n)}_{T_{ij}}(x)],\ \ \ \ x\in G.\]
For $\epsilon>0$, define
\[N_k(f,T)=\epsilon^{-1}\|f-\sigma_T^{(n)}\|_{\infty,k},\ \ \ f\in M_k(\mathcal{A}), T\in M_k(\mathcal{B}^{(n)}),\]
and denote $\mathcal{N}=(N_k)$.

\begin{lemma}\label{le1}
For any $T\in M_k(\mathcal{B}^{(n)})$, we have
\[L_k\left(\sigma_T^{(n)}\right)<L_k^{(n)}(T)+\epsilon.\]
\end{lemma}

\begin{proof}
Since $\sigma^{(n)}$ is a unital positive mapping from $\mathcal{B}^{(n)}$ to $\mathcal{A}$, $\sigma^{(n)}$ is unital completely positive and hence 
$\|\sigma^{(n)}\|_{cb}=1$ by Theorem 3.8 and Proposition 3.5 in \cite{pau}. So we have
\[\begin{array}{rcl}
&&L_k(\sigma_T^{(n)})\\
&=&\sup\left\{\frac{\left\|\left[{\displaystyle\lambda}_x\left({\displaystyle\sigma}^{(n)}_{T_{ij}}\right)-{\displaystyle\sigma}^{(n)}_{T_{ij}}\right]
\right\|_{\infty,k}}{{\displaystyle l(x)}}:
\ x\neq e, x\in G\right\}\\
&=&\sup\left\{\frac{\left\|\left[{\displaystyle\sigma}^{(n)}_{\left(\alpha_x^{(n)}(T_{ij})-T_{ij}\right)}\right]\right\|_{\infty,k}}{{\displaystyle l(x)}}:
\ x\neq e, x\in G\right\}\\
&\le&\sup\left\{\frac{\left\|\left[{\displaystyle\alpha}_x^{(n)}(T_{ij})-T_{ij}\right]\right\|}{{\displaystyle l(x)}}:
\ x\neq e, x\in G\right\}\\
&=&L_k^{(n)}(T)\\
&<&L_k^{(n)}(T)+\epsilon,\end{array}\]
by the $\alpha^{(n)}$-$\lambda$-equivariation of $\sigma^{(n)}$.
\end{proof}

Put on $\mathcal{A}$ the inner product from $L^2(G/H)$, while on $\mathcal{B}^{(n)}$ its Hilbert-Schmidt inner product. Then the mapping $\sigma^{(n)}$ from 
$\mathcal{B}^{(n)}$ to $\mathcal{A}$ has an adjoint operator $\breve{\sigma}^{(n)}$ from $\mathcal{A}$ to $\mathcal{B}^{(n)}$. For any $T\in\mathcal{B}^{(n)}$, 
a function $f\in \mathcal{A}$ such that $\breve{\sigma}^{(n)}_f=T$ is called a Berezin contravariant symbol for $T$. Moreover, $\breve{\sigma}^{(n)}$ is unital, positive, norm-nonincreasing, 
and $\lambda$-$\alpha^{(n)}$-equivariant (see page 75 in \cite{ri5}). From Theorem 3.10 and Proposition 3.5 in \cite{pau}, 
$\breve{\sigma}^{(n)}$ is unital completely positive and $\|\breve{\sigma}^{(n)}\|_{cb}=1$. 
So by the same argument as in the proof of Lemma \ref{le1}, we obtain:

\begin{lemma}\label{le2}
For any $f=[f_{ij}]\in M_k(\mathcal{A})$, we have
\[L_k^{(n)}\left(\breve{\sigma}^{(n)}_f\right)<L_k(f)+\epsilon,\]
where $\breve{\sigma}^{(n)}_f=[\breve{\sigma}^{(n)}_{f_{ij}}]\in M_k(\mathcal{B}^{(n)})$.
\end{lemma}

Denote
\[D_{L_k}(\varphi,\psi)=\sup\{\|\ll f,\varphi\gg-\ll f,\psi\gg\|: L_r(f)\le 1, f\in M_r(\mathcal{A}),r\in\mathbb{N}\},\]
for $\varphi,\psi\in CS_k(\mathcal{A})$, $k\in\mathbb{N}$, and 
\[h_{P^{(n)}}(x)=d^{(n)}\tau^{(n)}\left(P^{(n)}\alpha^{(n)}_x(P^{(n)})\right),\ \ x\in G/H,\]
where $\tau^{(n)}$ denotes the usual (un-normalized) trace on $\mathcal{B}^{(n)}$ and $d^{(n)}$ is the dimension of $\mathcal{H}^{(n)}$. Set 
\[\gamma^{(n)}=\int_{G/H}D_{L_1}(\hat{e},\hat{y})h_{P^{(n)}}(y)dy,\]
where every $y\in G/H$ is naturally identified with an element $\hat{y}$ of $CS_1(\mathcal{A})$. Then:

\begin{lemma}\label{le3}
For $f\in M_k(\mathcal{A})$, we have
\[\|f-\sigma^{(n)}(\breve{\sigma}^{(n)}_f)\|_{\infty,k}\le\gamma^{(n)}L_k(f).\]
\end{lemma}

\begin{proof}
Suppose $f=[f_{ij}]$. Then for any $x\in G/H$, we have
\[\begin{array}{rcl}
&&\|f(x)-(\sigma^{(n)}(\breve{\sigma}^{(n)}_f))(x)\|\\
&=&\left\|\left[\int_{G/H}\left(f_{ij}(x)-f_{ij}(y)\right)h_{P^{(n)}}(y^{-1}x)dy\right]\right\|\\
&=&\left\|\int_{G/H}\left[f_{ij}(x)-f_{ij}(y)\right]h_{P^{(n)}}(y^{-1}x)dy\right\|\\
&\le&\int_{G/H}\left\|\left[f_{ij}(x)-f_{ij}(y)\right]\right\|h_{P^{(n)}}(y^{-1}x)dy\\
&\le&L_k(f)\int_{G/H}D_{L_1}(\hat{x},\hat{y})h_{P^{(n)}}(y^{-1}x)dy\\
&=&L_k(f)\int_{G/H}D_{L_1}(\hat{e},\hat{y})h_{P^{(n)}}(y)dy\\
&\le&\gamma^{(n)}L_k(f),
\end{array}\]
by the formula (2.2) in \cite{ri5}. So 
\[\|f-\sigma^{(n)}(\breve{\sigma}^{(n)}_f)\|_{\infty,k}=\max\{\|f(x)-(\sigma^{(n)}(\breve{\sigma}^{(n)}_f))(x)\|: x\in G/H\}\le\gamma^{(n)}L_k(f).\]
\end{proof}

Since the sequence $\{\gamma^{(n)}\}$ converges to $0$ as $n\to\infty$ (see page 80 in \cite{ri5}), there is an $N_1\in\mathbb{N}$ such that 
$\gamma^{(n)}<\frac{\epsilon}2$ for $n> N^{(n)}_1$. So we obtain:

\begin{proposition}\label{pro1}
For $n> N_1$, $\mathcal{N}$ is a matrix bridge between $(\mathcal{B}^{(n)}, \mathcal{L}^{(n)})$ and $(\mathcal{A}, \mathcal{L})$, and hence 
$\mathcal{Q}=(Q_k)\in\mathcal{M}(\mathcal{L}^{(n)}, \mathcal{L})$, where
\[Q_k(f,T)=\max\{L_k^{(n)}(T), L_k(f), N_k(f,T)\}, \ \ \ (f,T)\in M_r(\mathcal{B}^{(n)}\oplus\mathcal{A}).\]
\end{proposition}

From Theorem 6.1 in \cite{ri5}, we have:

\begin{lemma}\label{le4}
There is an $N_2\in\mathbb{N}$ such that
\[\left\|T-\breve{\sigma}^{(n)}(\sigma^{(n)}_T)\right\|<\frac{\epsilon}2L^{(n)}_1(T),\]
for all $T\in\mathcal{B}^{(n)}$ and $n>N_2$.
\end{lemma}

\begin{theorem}
With notation as above, the quantized metric spaces $(\mathcal{B}^{(n)},\mathcal{L}^{(n)})$ converge to $(\mathcal{A},\mathcal{L})$ for 
quantized Gromov-Hausdorff distance as $n$ goes to $\infty$.
\end{theorem}

\begin{proof}
Given $\epsilon>0$. Choose $N=\max\{N_1,N_2\}$. Then for $n>N$, we have that $\mathcal{Q}\in\mathcal{M}(\mathcal{L}^{(n)}, \mathcal{L})$ by 
Proposition \ref{pro1}. Given $\varphi\in CS_k(\mathcal{A})$. we have $\varphi\circ\sigma^{(n)}\in CS_k(\mathcal{B}^{(n)})$, and
\[\begin{array}{rcl}
&&D_{L_k}(\varphi, \varphi\circ\sigma^{(n)})\\
&=&\sup\{\|\ll\varphi, f\gg-\ll\varphi\circ\sigma^{(n)}, T\gg\|: L_r(f,T)\le 1, \\
&&(f,T)\in M_r(\mathcal{A}\oplus\mathcal{B}^{(n)}), r\in\mathbb{N}\}\\
&=&\sup\{\|\ll\varphi, f-\sigma^{(n)}_T\gg\|: L_r(f,T)\le 1, 
(f,T)\in M_r(\mathcal{A}\oplus\mathcal{B}^{(n)}), r\in\mathbb{N}\}\\
&\le&\sup\{\|f-\sigma^{(n)}_T\|_{\infty,r}: L_r(f,T)\le 1, 
(f,T)\in M_r(\mathcal{A}\oplus\mathcal{B}^{(n)}), r\in\mathbb{N}\}\\
&\le&\epsilon.
\end{array}\]
On the other hand, if $\psi\in CS_k(\mathcal{B}^{(n)})$, then $\psi\circ\breve{\sigma}^{(n)}\in CS_k(\mathcal{A})$, and
\[\begin{array}{rcl}
&&D_{L_k}(\psi\circ\breve{\sigma}^{(n)},\psi)\\
&=&\sup\{\|\ll\psi\circ\breve{\sigma}^{(n)},f\gg-\ll\psi, T\gg\|: L_k(f,T)\le 1,\\
&&(f,T)\in M_k(\mathcal{A}\oplus\mathcal{B}^{(n)})\}\\
&=&\sup\{\|\ll\psi, \breve{\sigma}^{(n)}_f-T\gg\|: L_k(f,T)\le 1,
(f,T)\in M_k(\mathcal{A}\oplus\mathcal{B}^{(n)})\}\\
&\le&\sup\{\|\breve{\sigma}_f^{(n)}-T\|: L_k(f,T)\le 1, 
(f,T)\in M_k(\mathcal{A}\oplus\mathcal{B}^{(n)})\}\\
&\le&\sup\{\|\breve{\sigma}_f^{(n)}-\breve{\sigma}^{(n)}(\sigma^{(n)}_T)\|+\|\breve{\sigma}^{(n)}(\sigma^{(n)}_T)-T\|: L_k(f,T)\le 1, \\
&&(f,T)\in M_k(\mathcal{A}\oplus\mathcal{B}^{(n)})\}\\
&\le&\|f-\sigma^{(n)}_T\|_{\infty,k}+\sup\{\|\breve{\sigma}^{(n)}(\sigma^{(n)}_T)-T\|: L_k(f,T)\le 1,\\
&&(f,T)\in M_k(\mathcal{A}\oplus\mathcal{B}^{(n)})\}\\
&\le&\frac{\epsilon}2+\frac12 k^2\epsilon\\
&\le&k^2\epsilon,
\end{array}\]
by Lemma \ref{le4}. Therefore, for $n>N$, we have
\[\mathrm{dist}_{NC}(\mathcal{B}^{(n)},\mathcal{A})\le\epsilon,\]
that is, $\lim_{n\to\infty}\mathrm{dist}_{NC}(\mathcal{B}^{(n)},\mathcal{A})=0$.
\end{proof}

\section*{Acknowledgements}

I am grateful to Marc Rieffel for valuable discussions. I also thank Hanfeng Li for helpful comments. 
This research was partially supported by Shanghai Priority Academic Discipline, China Scholarship Council, and 
National Natural Science Foundation of China.

\bibliographystyle{amsplain}

\end{document}